\documentclass[11pt]{amsart}
\usepackage{geometry}                
\usepackage{fullpage}                   
\usepackage{graphicx}
\usepackage{amssymb}
\usepackage{epstopdf}
\usepackage{xcolor}
\usepackage[all,cmtip]{xy}
\newcommand{\Add}{\mathrm{Add}}

\newcommand{\PE}{\mathrm{PE}}
\newcommand{\Sim}{\mathrm{Sim}}
\newcommand{\pp}{\mathrm{pp}}

\newcommand{\M}{\mbox{\rm Mod-$R$}}

\DeclareMathOperator{\End}{End}

\DeclareMathOperator{\Tor}{Tor}
\DeclareMathOperator{\Img}{Im}
\DeclareMathOperator{\cY}{c\mathcal{Y}}
\newcommand{\im}{\mathrm{Im}\,}

\DeclareMathOperator{\soc}{soc}

\newcommand{\Ass}[1]{\mathrm{Ass} \,#1}
\newcommand{\Spec}[1]{\mathrm{Spec}(#1)}
\newcommand{\Max}[1]{\mathrm{Max}(#1)}

\newcommand{\C}{\ensuremath{\mathcal{C}}}

\newcommand{\p}{\ensuremath{\mathbf{p}}}

\newcommand{\q}{\ensuremath{\mathbf{q}}}
\newcommand{\tube}{\ensuremath{\mathbf{t}}}

\newcommand{\tensor}{\otimes}
\newcommand{\isom}{\cong}

\newcommand{\lra}{\longrightarrow}

\newcommand{\Supp}[1]{\mbox{\rm{Supp}} \,#1}

\newcommand{\Mod}[1]{\mathrm{Mod}\text{-}#1}
\newcommand{\LMod}[1]{#1\text{-}\mathrm{Mod}}
\renewcommand{\mod}[1]{\mathrm{mod}\text{-}#1}
\newcommand{\Lmod}[1]{#1\text{-}\mathrm{mod}}
\newcommand{\Ab}{\mathrm{Ab}}

\newcommand{\Der}[1]{\mathrm{D}(\Mod{#1})}
\newcommand{\D}[1]{\mathrm{D}(#1)}

\newcommand{\Hom}{\mathrm{Hom}}
\newcommand{\Ext}{\mathrm{Ext}}

\newcommand{\Ker}{\mathrm{Ker}\,}
\newcommand{\Coker}{\mathrm{Coker}}
\renewcommand{\Im}{\mathrm{Im}}

\newcommand{\id}[1]{\mathrm{Injdim}_R #1}

\newcommand{\Prod}[1]{\mathrm{Prod} #1}
\newcommand{\Cogen}[1]{\mathrm{Cogen} #1}

\newcommand{\Inj}{\mathrm{Inj}}
\newcommand{\Pinj}{\mathrm{Pinj}}

\newcommand{\Acal}{\mathcal{A}}
\newcommand{\Bcal}{\mathcal{B}}
\newcommand{\Ccal}{\mathcal{C}}
\newcommand{\Dcal}{\mathcal{D}}
\newcommand{\Ecal}{\mathcal{E}}
\newcommand{\Fcal}{\mathcal{F}}
\newcommand{\Gcal}{\mathcal{G}}
\newcommand{\Hcal}{\mathcal{H}}
\newcommand{\Ical}{\mathcal{I}}

\newcommand{\Kcal}{\mathcal{K}}

\newcommand{\Mcal}{\mathcal{M}}
\newcommand{\Ncal}{\mathcal{N}}

\newcommand{\Qcal}{\mathcal{Q}}
\newcommand{\Rcal}{\mathcal{R}}
\newcommand{\Scal}{\mathcal{S}}
\newcommand{\Tcal}{\mathcal{T}}
\newcommand{\Ucal}{\mathcal{U}}
\newcommand{\Vcal}{\mathcal{V}}
\newcommand{\Wcal}{\mathcal{W}}
\newcommand{\Xcal}{\mathcal{X}}
\newcommand{\Ycal}{\mathcal{Y}}

\newtheorem{theorem}{Theorem}[section]
\newtheorem{definition}[theorem]{Definition}
\newtheorem{lemma}[theorem]{Lemma}
\newtheorem{corollary}[theorem]{Corollary}                
\newtheorem{proposition}[theorem]{Proposition}
\newtheorem{remark}[theorem]{Remark}
\newtheorem{notation}[theorem]{Notation}
\newtheorem{example}[theorem]{Example}

\numberwithin{equation}{section}
\title{Simples in a cotilting heart}
\author{Lidia Angeleri H{\"u}gel, Ivo Herzog and Rosanna Laking}
\address{Lidia Angeleri H\"ugel, Dipartimento di Informatica - Settore di Matematica, Universit\`a degli Studi di Verona, Strada le Grazie 15 - Ca' Vignal, I-37134 Verona, Italy} 
\email{lidia.angeleri@univr.it}
\address{Ivo Herzog, 4240 Campus Drive, The Ohio State University, Lima, OH 45804, USA}
\email{herzog.23@osu.edu}
\address{Rosanna Laking, Dipartimento di Informatica - Settore di Matematica, Universit\`a degli Studi di Verona, Strada le Grazie 15 - Ca' Vignal, I-37134 Verona, Italy} 
\email{rosanna.laking@univr.it}
\thanks{Acknowledgments: The authors would like to thank the referee for valuable suggestions.  The first and third named authors acknowledge support from the project \textit{REDCOM: Reducing complexity in algebra, logic, combinatorics}, financed by the programme  \textit{Ricerca Scientifica di Eccellenza 2018} of the Fondazione Cariverona.  The third named author was also supported by the European Union's Horizon 2020 research and innovation programme under the Marie Sk{\l}odowska-Curie Grant Agreement No.~797281. The second named author was partially supported by \textit{NSF Grant DMS 12-01523.}}
\date{\today}                                           

\begin{document}
\maketitle

\begin{abstract}
Every cotilting module over a ring $R$ induces a t-structure  with a Grothendieck heart in the derived category $\Dcal(\Mod R)$. We determine the simple objects in this heart and their injective envelopes, combining torsion-theoretic aspects with methods from the model theory of modules and Auslander-Reiten theory.
\end{abstract} 

\section{Introduction}
The notion of a $t${\em -structure} $\tau$ on a triangulated category $\Tcal$ appears in the work of Beilinson, Bernstein, and Deligne~\cite{beilinson:bernstein:deligne:1981} as a means to associate to $\Tcal$ an abelian category, which then arises as the {\em heart} $\Hcal_{\tau}$ of the $t$-structure. For example, if the triangulated category is the derived category $\Dcal (\Acal)$ of an abelian category $\Acal,$ an appropriately chosen $t$-structure $\tau$ will recover, according to this process, the given abelian category $\Hcal_{\tau} \isom \Acal.$ Other choices of $t$-structure on $\Dcal (\Acal)$ will give rise to hearts that may be derived equivalent to $\Acal.$ 

The primary aim of~\cite{beilinson:bernstein:deligne:1981} was to introduce the abelian category of {\em perverse sheaves} (see~\cite[Ch 8]{kirwan:woolf:2006}) on a topological pseudomanifold $X$ of even dimension, equipped with a stratification with no odd-dimensional strata. The triangulated category was the bounded derived category of sheaves on $X$ and the $t$-structure was chosen to yield the perverse sheaves on $X$ as the objects of the heart $\Hcal_{\tau}.$ As these perverse sheaves were seen to have finite length, attention naturally turned to the {\em simple} ones, which were determined~\cite[Theorem 4.3.1]{beilinson:bernstein:deligne:1981} by the intersection homology of the connected strata, appropriately shifted (see also~\cite[Theorem 8.1.8]{kirwan:woolf:2006}).

In the same spirit, our interest in this paper are the simple objects of the heart of the {\em HRS tilt} of a torsion pair $(\Tcal, \Fcal)$ in the category $\Mod{R}$ of modules over a ring. The HRS tilt of $(\Tcal, \Fcal)$ is a $t$-structure $\tau$ introduced by 
Happel, Reiten, and Smal\o~\cite{happel:reiten:smaloe:1996} on the derived category $\Dcal (\Mod{R}).$ The objects of the heart $\Hcal_{\tau}$ need not all be of finite length, but Happel, Reiten and Smal{\o} showed that $\Hcal_{\tau}$ contains a torsion pair $(\Fcal, \Tcal[-1])$ whose constituent classes are equivalent to those of $(\Tcal, \Fcal),$ with the roles reversed. When viewed with regard to the torsion pair $(\Fcal, \Tcal [-1]),$ every simple object of $\Hcal_{\tau}$ is evidently either torsion or torsionfree. 

We provide a {\em torsion theoretic} description of the simple objects of $\Hcal_{\tau}$ using the notion of an \textbf{almost torsionfree} module (Definition~\ref{def: tf/t t/tf}) and its dual, that of an \textbf{almost torsion} module. Every torsionfree module is almost torsionfree, but there may be others, which are necessarily torsion. The dual statement also holds and we characterise in Theorem~\ref{thm: simples t/tf tf/t} the torsionfree simple objects of $\Hcal_{\tau}$ to be those of the form $T[-1]$ where $T_R \in \Tcal \subseteq \Mod{R}$ is an almost torsionfree torsion module, and the torsion simples of the heart to be the objects that correspond to almost torsion torsionfree modules.\bigskip

If $C \in \Mod{R}$ is a $1$-cotilting module, then the cotilting class $\Ccal = {^{\perp_1}}C = \Cogen (C)$ is a torsionfree class in $\Mod{R}$ and we call the heart of the HRS tilt of the torsion pair $(\Qcal, \Ccal)$ a \textbf{cotilting heart.} Every cotilting heart is a Grothendieck category~\cite{colpi:gregorio:mantese:2007} whose subcategory $\Inj (\Hcal_{\tau})$ of injective objects is known to be equivalent to $\Prod (C) \subseteq \Ccal.$ As such, the injective objects of $\Hcal_{\tau}$ are torsion, but when $\Prod (C) \subseteq \Ccal \subseteq \Mod{R}$ is regarded as a subcategory of $R$-modules, it consists of torsionfree modules. As cotilting modules are pure-injective, so are all the objects of $\Prod (C).$ Corollary~\ref{cor: inj env neg-iso} makes use of the notion of a \textbf{neg-isolated} indecomposable pure-injective module (\S \ref{subsec: inj env of simples}) from the {\em model theory of modules} to characterise the injective envelopes of simple objects of the heart, when they are considered as modules in the definable subcategory $\Ccal \subseteq \Mod{R}.$ It states that they are precisely the neg-isolated indecomposable pure-injectives of $\Ccal$ that belong $\Prod (C).$

Among the neg-isolated indecomposable pure-injective modules of a definable subcategory of $\Mod{R}$ such as $\Ccal,$ there is the distinguished class of \textbf{critical} neg-isolated indecomposable pure-injectives $U,$ determined by the property that every monomorphism $U \to V$ in $\Mod{R}$ with $V \in \Ccal$ is a split monomorphism (Proposition  \ref{prop: critical properties}).  These are the torsionfree modules that correspond to the injective envelopes of torsion simple objects of the heart. It is a general fact 
about definable subcategories that there exist enough critical neg-isolated indecomposables, in the sense that every torsionfree module $F \in \Ccal$ may be embedded - not necessarily purely - into a direct product of critical neg-isolated indecomposables in $\Ccal.$ It follows that no $1$-cotilting module is superdecomposable. 

The characterisations of the torsionfree and torsion simple objects of a cotiliting heart in terms of almost torsionfree and almost torsion modules are categorically dual and seem to give the two kinds of simple object equal status. 
The question of existence however does not. We call the neg-isolated indecomposable pure-injectives of $\Ccal$ that correspond to injective envelopes of torsionfree simple objects of the heart \textbf{special.} In stark contrast to the critical neg-isolated indecomposables, there is a $1$-cotilting module $C_{\Lambda}$ over the Kronecker algebra $\Lambda$ (Example~\ref{no t/tf}) whose cotilting class contains no special neg-isolated indecomposable pure-injectives. In other words, every simple object of the cotilting heart of $C_{\Lambda}$ is torsion.
\bigskip

All of our characterisations of the simple objects of a cotilting heart may be regarded as part of {\em Auslander-Reiten theory,} but only the last makes direct appeal to almost split morphisms in the module category. This final description relies on the {\em approximation theory} of the complete cotorsion pair $(\Ccal, \Ccal^{\perp}).$ The almost split morphisms that appear are left almost split morphisms that enjoy the \textbf{strong} uniqueness property (Definition~\ref{D:strong left AR seq}). It is included in the following summary of all our results on the torsion simple objects of a cotilting heart.
\bigskip

\noindent {\bf Theorem A} (Theorems~\ref{thm: simples t/tf tf/t} and \ref{Thm: las and specials}, Proposition~\ref{prop: inj envelope}, Corollary~\ref{cor: crit inj env}) 

{\it The following statements are equivalent for a module $N$.
\begin{enumerate}
\item $N$ is isomorphic to the injective envelope of a torsion simple $S$ in $\Hcal_\tau$.
\item $N$ is a critical neg-isolated module in $\Ccal$.
\item There exists a short exact sequence $$ \xymatrix@1{ 0 \ar[r]  & S   \ar[r]^a  & N  \ar[r]^b  & \bar{N}   \ar[r]  & 0 }$$
in $\Mod{R},$ where $S$ is torsionfree, almost torsion, $a$ is a $\Ccal^{\perp_1}$-envelope, and $b$ is a strong left almost split morphism in $\Ccal$.
\end{enumerate}
} 

\bigskip

\noindent A strong left almost split morphism is either a monomorphism or epimorphism (Lemma~\ref{lem: strong las mono or epi}). Theorem A includes a characterisation of the torsion simple objects of the heart as the torsionfree modules that appear as {\em kernels} of strong left almost split morphisms in $\Ccal,$ while  its dual, the next Theorem B, characterises the torsionfree simple objects as shifts of torsion modules that arise as {\em cokernels} of strong left almost split morphisms in $\Ccal.$
\bigskip

{\bf Theorem B} (Theorems~\ref{thm: simples t/tf tf/t} and \ref{Thm: las and specials}, Proposition~\ref{prop: inj envelope}, Proposition~\ref{prop: special}) 

{\it The following statements are equivalent for a module $N$.
\begin{enumerate}
\item $N$ is isomorphic to the injective envelope of a torsionfree simple $S[-1]$ in $\Hcal_\tau$.
\item $N$ is a special neg-isolated module in $\Ccal$.
\item There exists a short exact sequence $$ \xymatrix@1{ 0 \ar[r]  & N   \ar[r]^a  & \bar{N}  \ar[r]^b  & S   \ar[r]  & 0 }$$ 
in $\Mod{R},$ where $S$ is torsion, almost torsionfree, $a$ is a a strong left almost split morphism in $\Ccal$, and  $b$ is a  $\Ccal$-cover.
\end{enumerate}
}

\bigskip

The simple objects in cotilting hearts are crucial to understanding the phenomenon of mutation and to describe the lattice tors-$R$ of torsion classes in the category $\mod{R}$ of finite dimensional modules over a finite dimensional algebra $R$. Indeed, the  simple objects in the heart $\Hcal_\tau$ correspond  to the arrows in the Hasse quiver of tors-$R$ which are incident to the torsion class $\Qcal\cap\mod{R}$, or equivalently, to the irreducible mutations of the cotilting module $C$, cf.~\cite{demonet:iyama:reading:reiten:thomas:2017,barnard:carroll:zhu:2019,angeleri:laking:stovicek:vitoria:2022}.  In a forthcoming paper \cite{angeleri:laking:sentieri:2022}, we will employ Theorems A and B to  obtain an explicit description of mutation of cotilting (or more generally, cosilting) modules. This will allow us to interpret mutation  as an  operation 
on the Ziegler spectrum of $R$ which will amount to replacing  critical neg-isolated summands by special ones, or viceversa. 

\bigskip

\section{Background}

\subsection{Notation}
In this section we fix our basic notations and conventions.

Let $R$ be a unital associative ring.  We denote the category of right $R$-modules by $\Mod{R}$ and the category of left $R$-modules by $\LMod{R}$.  The full subcategories of finitely presented modules are denoted $\mod{R}$ and $\Lmod{R}$ respectively.  The derived category of $\Mod{R}$ is denoted $\Der{R}$.  We abbreviate the Hom-spaces in $\Der{R}$ in the following way: \[\Hom_{\D{R}}(X, Y) := \Hom_{\Der{R}}(X,Y)\] for all complexes $X$, $Y$. 

All subcategories will be strict (i.e.~closed under isomorphisms) and, for a full subcategory $\Bcal$, we will use the notation $B \in \Bcal$ to indicate that $B$ is an object of $\Bcal$.

Let $\Xcal$ be a set of objects in an additive category $\Acal$ with products.  Then we use the notation $\Prod(\Xcal)$ for the set of direct summands of products of copies of objects contained in $\Xcal$.  In the case where $\Acal$ is Grothendieck abelian, we will use $\Cogen(\Xcal)$ to denote the set of subobjects of objects contained in $\Prod(\Xcal)$.  We will write $\Inj(\Xcal)$ for the class of injective objects in the category $\Acal$ that are contained in $\Xcal$.  We will consider the following full perpendicular subcategories determined by a subset  $I\subseteq\{0,1\}$: \[ \Xcal^{\perp_I} :=  \{ M \in  \Acal \mid \Ext^i_\Acal(X, M) = 0 \text{ for all } X\in \Xcal \text{ and } i\in I\}\] 
\[  {}^{\perp_I}\Xcal :=  \{ M \in  \Acal \mid \Ext^i_\Acal(M, X) = 0 \text{ for all } X\in \Xcal \text{ and } i\in I \}.\]  In the case where $\Xcal = \{X\}$, we will use the notation $X^{\perp_I}$ for $\Xcal^{\perp_I}$ and $\Prod(X)$ for $\Prod(\Xcal)$ etc. Furthermore, we will often just write $\Xcal^{\perp_0}$ instead of $\Xcal^{\perp_{\{0\}}}$ etc.

\subsection{Torsion pairs and HRS-tilts}

In this subsection we introduce the notion of an HRS-tilt, due to Happel, Reiten and Smal{\o}.  The idea of their work is to produce a t-structure in the derived category $\Der{R}$ from a given torsion pair in $\Mod{R}$.  More details about the construction and properties of this t-structure can be found in \cite{happel:reiten:smaloe:1996}.\newline

Torsion pairs, first introduced by Dickson \cite{dickson:1966}, will be a central object of study in the latter sections of this article.  The following is the definition of a torsion pair in abelian category $\Acal$.

\begin{definition}
A pair of full subcategories $(\Tcal, \Fcal)$ of $\Acal$ is called a \textbf{torsion pair} if the following conditions hold.\begin{enumerate}
\item For every $T \in \Tcal$ and $F \in \Fcal$, we have that $\Hom_\Acal(T, F)=0$.
\item For every $X$ in $\Acal$, there exists a short exact sequence \[ 0 \to t(X) \to X \to X/t(X) \to 0 \] where $t(X) \in \Tcal$ and $X/t(X) \in \Fcal$.
\end{enumerate}  
We call $\Tcal$ the \textbf{torsion class} and $\Fcal$ the \textbf{torsionfree class}.  If, in addition, the class $\Tcal$ is closed under subobjects, then the torsion pair is called \textbf{hereditary}.
\end{definition}

We extend the above terminology to objects: the objects $T$ in $\Tcal$ are called \textbf{torsion} and the objects $F$ in $\Fcal$ are called \textbf{torsionfree}.

The next result shows that such a torsion pair in $\Mod{R}$ yields a t-structure in $\Der{R}$, in the sense of \cite{beilinson:bernstein:deligne:1981}.  Note that we define our t-structure to consist of two Hom-orthogonal classes; this differs from the original definition by a shift. 

\begin{proposition}[{\cite[Prop.~I.2.1]{happel:reiten:smaloe:1996}}]
Let $\tau = (\Tcal, \Fcal)$ be a torsion pair in $\Mod{R}$.  The two full subcategories \[ \Ucal_\tau = \{ X \in \Der{R} \mid H^0(X) \in \Tcal, H^i(X)=0 \text{ for } i>0\}\] \[ \Vcal_\tau = \{X\in \Der{R} \mid H^{0}(X) \in \Fcal, H^i(X)=0 \text{ for } i<0\}\] of $\Der{R}$ form a t-structure.
\end{proposition}

 We will refer to this t-structure as the \textbf{HRS-tilt of $(\Tcal, \Fcal)$}.  It is shown in \cite{beilinson:bernstein:deligne:1981} that the \textbf{heart} $\Hcal_\tau := \Ucal_\tau[-1] \cap \Vcal_\tau$ of the t-structure $(\Ucal_\tau, \Vcal_\tau)$ is an abelian category whose short exact sequences $0 \to X \to Y \to Z \to 0$ are given by the triangles $X \to Y \to Z \to X[1]$ of $\Der{R}$ such that $X$, $Y$ and $Z$ are contained in $\Hcal_\tau$.  For any two objects $X$ and $Y$ in $\Hcal_\tau$, there are functorial isomorphisms \[ \Hom_{\D{R}}(X, Y[i]) \cong \Ext^i_{\Hcal_\tau}(X, Y) \text{ for } i=0,1.\]

\noindent Moreover, $(\Fcal,\Tcal[-1])$ is a torsion pair in $\Hcal_\tau$ by \cite[Cor.~I.2.2]{happel:reiten:smaloe:1996}.\newline

We will make use of the following lemma in Section \ref{sec: simples}.

\begin{lemma}\label{heartlemma}
Let $\tau = (\Tcal, \Fcal)$ be a torsion pair in $\Mod{R}$.
 \begin{enumerate}
\item  Let $f\colon X\to Y$ be a morphism in $\Hcal_\tau$, and let $Z$ be the cone of $f$ in $\Der{R}$. Consider the canonical triangle  \[K\to Z\to W\to K[1]\]
where $K\in \Ucal_\tau$ and  $W\in\Vcal_\tau.$
Then \[\Ker_{\Hcal_\tau} (f)=K[-1],\quad \Coker_{\Hcal_\tau} (f)=W.\]
\item Let $h\colon Y\to X$ be an $R$-homomorphism with $Y,X\in\Fcal$. The morphism  $h$  is a monomorphism in $\Hcal_\tau$ if and only if $\Ker(h)=0$ and $\Coker(h)\in\Fcal$, and $h$ is an epimorphism in $\Hcal_\tau$ if and only if $\Coker(h)\in\Tcal$.

\item Let $h:Y\to X$ be a $R$-homomorphism with $Y,X\in\Tcal$. The morphism  $h[-1]$ is a monomorphism in $\Hcal_\tau$ if and only if $\Ker(h)\in\Fcal$, and $h[-1]$ is an epimorphism in $\Hcal_\tau$  if and only if $\Coker(h)=0$ and $\Ker(h)\in\Tcal$.
\end{enumerate}
\end{lemma}
\begin{proof}  Recall that the cone of a morphism $h$ in $\Mod{R}$ has  homologies $\Ker(h)$ in degree $-1$, $\Coker(h)$ in degree $0$, and zero elsewhere. 

(1) This is a standard property of t-structures.  See, for example, \cite[pp.281]{gelfand:manin:2003}.

(2) We know from (1)  that  $\Ker_{\Hcal_\tau}(h)=0$ if and only if the cone of $h$ belongs to $\Vcal_\tau$. This means $\Ker(h)=0$ and $\Coker(h)\in\Fcal$. Similarly, $\Coker_{\Hcal_\tau}(h)=0$ if and only if the cone of $h$ belongs to $\Ucal_\tau$, which means that $\Coker(h) \in\Tcal$. 

(3) The cone of $h[-1]$ belongs to $\Vcal_\tau$ if and only if  $\Ker(h)\in\Fcal$, and it belongs to  $\Ucal_\tau$  if and only if $\Coker(h)=0$ and $\Ker(h)\in\Tcal$.
\end{proof}

\subsection{Cotilting modules and cotorsion pairs}

In this paper we will focus on HRS-tilts of torsion pairs induced by cotilting modules.  We now introduce these modules and collect together some of their important properties.  The definition of a (possibly infinitely generated) cotilting module first appeared in \cite{colpi:d'este:tonolo:1997}, dualising the definition of \cite{colpi:trlifaj:1995}.

\begin{definition}\label{Def: cotilting}
A right $R$-module $C$ is called a \textbf{cotilting module} if the following three statements hold.\begin{enumerate}
\item $\id{C} \leq 1$.
\item $\Ext^1_R(C^\kappa, C) = 0$ for all cardinals $\kappa$.
\item There exists a short exact sequence $0 \to C_1 \to C_0 \to I \to 0$ where $C_i \in \Prod(C)$ for $i=0,1$ and $I$ is an injective cogenerator of $\Mod{R}$.
\end{enumerate}
We say that cotilting modules $C$ and $C'$ are \textbf{equivalent} if $\Prod(C) = \Prod(C')$.
\end{definition}

In \cite[Prop.~1.7]{colpi:d'este:tonolo:1997}, the authors show that $\Cogen(C) = {}^{\perp_1}C$ and, moreover, that this equality characterises cotilting modules. We call this class $\Ccal := \Cogen(C) = {}^{\perp_1}C$ the \textbf{cotilting class} associated to $C$ and it follows that $\tau = (\Qcal, \Ccal) := ({}^{\perp_0}C, \Cogen(C))$ is a (faithful) torsion pair.  We call the heart of the HRS-tilt of $\tau$ the associated \textbf{cotilting heart}. 

We know from \cite{colpi:gregorio:mantese:2007} that a cotilting heart $\Hcal_\tau$ is a Grothendieck category with injective cogenerator $C$ so, in particular, we have $\Inj(\Hcal_\tau) = \Prod({C})$. 

\begin{remark}
Often the term cotilting module is used for the more general notion of an $n$-cotilting module, which was first defined in \cite{angeleri:coelho:2001}.  In that context, the modules specified in Definition \ref{Def: cotilting} are called $1$-cotilting modules.  Since we will not be considering $n$-cotilting modules for $n>1$, we will use the term cotilting module to refer to a $1$-cotilting module.
\end{remark}

It was shown in \cite{bazzoni:2003} that every cotilting module is pure-injective and every cotilting class is definable (see Sections \ref{Sec: functor cats} and \ref{sec:pi adj} for definitions of these terms).  As a consequence, the class $\Ccal$ is closed under direct limits, and the cotorsion pair $(\Ccal, \Ccal^{\perp_1}) = ({}^{\perp_1}C, ({}^{\perp_1}C)^{\perp_1})$ cogenerated by $C$ is a perfect cotorsion pair.  In particular, for every module $M$ in $\Mod{R}$, there exist special approximation sequences \[ 0 \longrightarrow X \longrightarrow Y \overset{a}{\longrightarrow} M \longrightarrow 0 \] \[ 0 \longrightarrow M \overset{b}{\longrightarrow} X' \longrightarrow Y' \longrightarrow 0  \] such that $X, X' \in \Ccal^{\perp_1}$ and $Y, Y' \in \Ccal$.  In particular, $a$ is a $\Ccal$-cover and $b$ is a $\Ccal^{\perp_1}$-envelope.  Moreover, we have that $\Ccal \cap \Ccal^{\perp_1} = \Prod(C)$.  For more details on covers, envelopes and cotorsion pairs, we refer the reader to \cite{goebel:trlifaj:2006}.

\subsection{Injective envelopes of simples and left almost split morphisms}

In this section we will prove some preliminary results  connecting simple objects in a cotilting heart to left almost split morphisms. Our considerations are inspired by \cite{crawley-boevey:1992}.
\begin{definition} \label{D:strong left AR seq}
Let $\Xcal$ be an additive category.  A morphism $f\colon X \to Y$ in $\Xcal$ is called a \textbf{left almost split morphism} if it is not a split monomorphism and, for any $g \colon X \to Z$ that is not a split monomorphism, there exists a morphism $h \colon Y \to Z$ such that $g = hf$.  If the morphism $h$ is unique for every such $g$, then we call $f$ a \textbf{strong} left almost split morphism.
\end{definition}

We begin with the following general result about Grothendieck abelian categories.

\begin{proposition}\label{prop: las inj env simple}
Let $\Gcal$ be a Grothendieck abelian category and, for any object $M$ in $\Gcal$, let $E(M)$ denote the injective envelope of $M$. \begin{enumerate}
\item If $S$ is a simple object, then the canonical morphism $E(S) \to E(E(S)/S)$ is a left almost split morphism in $\Inj(\Gcal)$.
\item If $f \colon E \to E^+$ is a left almost split morphism in $\Inj(\Gcal)$, then the kernel $\Ker(f)$ is simple and the canonical embedding $\Ker(f) \to E$ is the injective envelope of $\Ker(f)$.
\item If $f \colon E \to E^+$ is a left almost split morphism in $\Inj(\Gcal)$, then the canonical epimorphism $g: E \to \Im(f)$ is a strong left almost split morphism in $\Gcal$.
\item If $g \colon E \to \tilde{E}$ is a left almost split morphism in $\Gcal$ with $E$ in $\Inj(\Gcal)$, and $e \colon \tilde{E} \to E(\tilde{E})$ is the injective envelope of $\tilde{E}$, then $f:=eg$ is a left almost split morphism in $\Inj(\Gcal)$.
\end{enumerate}
\end{proposition}
\begin{proof}
\textbf{(1)}   Consider the morphism $f$ given by the composition of the quotient $E(S) \to E(S)/S$ with the injective envelope $E(S)/S \to E(E(S)/S)$.  This morphism is not a split monomorphism.  Any other morphism $g \colon E(S)\to F$ in $\Inj(\Gcal)$ that is not a split monomorphism must have a non-trivial kernel $K$ and so $K$ necessarily contains $S$ because $S$ is essential in $E(S)$.  It follows that $g$ factors through $f$ as required.\smallskip

\textbf{(2)} Consider the kernel $0 \to K \overset{k}{\to} E \overset{f}{\rightarrow} E^+$ of $f$ in $\Gcal$.  We will show that $K = \Ker(f)$ is simple.  Clearly $K \neq 0$ because $f$ is not a split monomorphism. 
Moreover,  every non-zero subobject $G \subset K$ coincides with $K$, because the composition of the quotient $E \to E/G$ with the injective envelope $E/G \to E(E/G)$ of $E/G$ is not a split monomorphism and thus factors through $f$.

Let $e \colon K \to E(K)$ be the injective envelope of $K$.  Since $k \colon K \to E$ is a monomorphism and $E$ is injective, there exists a split epimorphism $m \colon E \to E(K)$ such that $e = mk$.  If $m$ is not a monomorphism, then there exists a morphism $g \colon E^+ \to E(K)$ such that $gf = m$.  This implies that $0 = gfk = mk = e$, which is a contradiction.  Therefore $m$ is an isomorphism.  \smallskip

\textbf{(3)} By (2), we have an exact sequence $0\to S \overset{i}{\to} E \overset{f}{\to} E^+$ where $i$ is the injective envelope of $S$ and $S$ is simple.  Consider the short exact sequence \[ 0 \to S \overset{i}{\to} E \overset{g}{\to} E/S \to 0.\]  We will show that $g$ is a strong left almost split morphism in $\Gcal$.  Note that $g$ is not a monomorphism and so cannot be a split monomorphism.  Consider a morphism $a \colon E \to M$ that is not a split monomorphism.  If $ai \neq 0$, then $ai$ must be a monomorphism because $S$ is simple. Then $a$ is a monomorphism because $i$ is an essential monomorphism.  This implies that $a$ splits because $E$ is injective, but this is a contradiction.  Thus $ai =0$ and therefore $a$ factors uniquely through the cokernel $E/S \cong \Im(f)$ of $i$, as required.\smallskip

\textbf{(4)} Consider $f := eg$ where $e\colon \tilde{E} \to E({\tilde{E}})$ is the injective envelope of $\tilde{E}$.  We will show that $f$ is a left almost split morphism in $\Inj(\Gcal)$.  Firstly, $f$ is not a split monomorphism because otherwise $g$ is a monomorphism and therefore split (since $E$ is injective).  Let $a \colon E \to E'$ be a morphism in $\Inj(\Gcal)$ that is not a split monomorphism.  As $g$ is a left almost split morphism in $\Gcal$, we have that there exists a morphism $b \colon \tilde{E} \to E'$ such that $a = bg$.  Moreover, since $E'$ is injective and $e$ is a monomorphism, we have that there exists a morphism $c \colon E(\tilde{E}) \to E'$ such that $a = c(eg) = cf$, as required.
\end{proof}

\begin{remark}\label{rem: las in cotilting class}
Following all the notation of Proposition \ref{prop: las inj env simple}, assume that the Grothendieck category $\Gcal = \Hcal_\tau$ is a cotilting heart with respect to the cotilting torsion pair $\tau = (\Qcal, \Ccal)$. 
Then, in the argument for Proposition \ref{prop: las inj env simple}(3), the object $\Im_{\Hcal_\tau}(f)$ is in $\Ccal$ because $\Ccal$ is a torsion class. Hence 
$g$ is a strong left almost split morphism in the subcategory $\Ccal$. Moreover,  the argument for Proposition \ref{prop: las inj env simple}(4) only requires that $g$ is a  left almost split morphism in $\Ccal$ since the injective objects in $\Inj(\Hcal_\tau) = \Prod(C)$ are contained in $\Ccal$. So, every left almost split morphism $E\to \tilde{E}$ in $\Ccal$ with $E\in\Prod(C)$ induces a left almost split morphism in $\Inj(\Hcal_\tau)$.
\end{remark}

\begin{corollary}\label{cor: inj env simp}
Let $\Gcal$ be a Grothendieck abelian category and let $\Inj(\Gcal)$ denote the full subcategory of injective objects in $\Gcal$.  The following statements are equivalent for an object $E$ of $\Inj(\Gcal)$. \begin{enumerate}
\item $E$ is isomorphic to the injective envelope $E(S)$ of a simple object $S$ in $\Gcal$.
\item There exists a left almost split morphism $f \colon E \to E^+$ in $\Inj(\Gcal)$.
\item There exists a strong left almost split morphism $g\colon E \to \tilde{E}$ in $\Gcal$.
\item There exists a left almost split morphism $g\colon E \to \tilde{E}$ in $\Gcal$.
\end{enumerate}  
\end{corollary}

Proposition \ref{prop: las inj env simple} and Remark \ref{rem: las in cotilting class} yield the following corollary in the special case where $\Gcal$ is a cotilting heart.

\begin{corollary}\label{cor: las simples}
Let $\tau = (\Qcal, \Ccal)$ be a cotilting torsion pair in $\Mod{R}$ and let $\Inj(\Hcal_\tau)$ denote the full subcategory of injective objects in $\Hcal_\tau$.  The following statements are equivalent for an object $E$ of $\Inj(\Hcal_\tau)$. \begin{enumerate}
\item $E$ is isomorphic to the injective envelope $E(S)$ of a simple object $S$ in $\Hcal_\tau$.
\item There exists a left almost split morphism $f \colon E \to E^+$ in $\Inj(\Hcal_\tau)$.
\item There exists a (strong) left almost split morphism $g\colon E \to \tilde{E}$ in $\Hcal_\tau$.
\item There exists a (strong) left almost split morphism $g\colon E \to \bar{E}$ in the torsion class $\Ccal$.
\end{enumerate}  
\end{corollary}

\subsection{Localisation in abelian Grothendieck categories} \label{sec:Gr loc}
A torsion pair $(\Tcal, \Fcal)$ in an abelian category is \textbf{hereditary} if the torsion class $\Tcal$ is closed under subobjects. If the abelian category is Grothendieck, this is equivalent to the torsionfree class being closed under injective envelopes. 
In that case, the torsion pair $(\Tcal, \Fcal)$ is cogenerated by its torsionfree injective objects. Define the \textbf{localisation} $\Gcal/ \Tcal$ of $\Gcal$ at $\Tcal$ to be the category whose objects are the same as the objects $X$ of $\Gcal,$ 
but denoted by $X_{\Tcal}.$ The morphisms between two objects are given by the set $\Hom_{\Gcal/\Tcal} (X_{\Tcal},Y_{\Tcal}) := \varinjlim \; \Hom_{\Gcal} (X', Y/Y')$ where $X'$ ranges over the subobjects of $X$ such that $X/X' \in \Tcal$ and $Y'$ ranges over the subobjects of $Y$ such that $Y'\in \Tcal$. The work of Gabriel shows that the \textbf{localisation functor} 
$L_{\Tcal} \colon \Gcal \to \Gcal/\Tcal,$ $X \mapsto X_{\Tcal},$ is the left adjoint of an adjunction
\begin{equation} \xymatrix@C=15pt{\Gcal \ar@/^1pc/[rr]^{L_{\Tcal}} & \bot &  \Gcal/\Tcal \ar@/^1pc/[ll]^{R_{\Tcal}}.} \end{equation}
The adjoint property allows us to calculate hom groups in the localisation: if $X \in \Gcal$ and $Y_\Tcal \in \Gcal/\Tcal,$ then $\Hom_{\Gcal/\Tcal} (X_{\Tcal}, Y_{\Tcal}) \isom \Hom_{\Gcal} (X, R_{\Tcal}(Y_{\Tcal})).$ 

The {\em left} adjoint $L_{\Tcal}$ is exact and the {\em right} adjoint $R_{\Tcal} \colon \Gcal/\Tcal \to \Gcal$ is fully faithful. We may therefore identify the localisation category $\Gcal/\Tcal$ with the full subcategory of $\Gcal$ given by the image of 
$R_{\Tcal}$.  We note that $\Gcal/\Tcal$ is contained in $\Fcal$. Because the right adjoint of an exact functor preserves injective objects, we may regard $\Inj (\Gcal/\Tcal)$ under this identification as a subcategory of $\Inj (\Gcal);$ it is precisely the subcategory 
$\Inj (\Gcal/\Tcal) = \Inj (\Gcal)\, \cap \, \Fcal$ of torsionfree injective objects. 

As the right adjoint $R_{\Tcal}$ is left exact, we may identify the entire localisation $\Gcal/\Tcal$ with the equivalent subcategory $\Cogen^2 (\Inj (\Gcal/\Tcal)) \subseteq \Gcal$ consisting of the objects in $\Gcal$ with a copresentation by torsionfree injectives. 
For more details on localisation in Grothendieck categories, we refer the reader to \cite[Ch.~4]{popescu:1973}.

\section{Simple objects in the heart}\label{sec: simples}

In this section we consider the simple objects in the heart $\Hcal_\tau$ of the HRS-tilt of a torsion pair $\tau = (\Tcal, \Fcal)$ in $\Mod{R}$.  Since $(\Fcal, \Tcal[-1])$ is a torsion pair in $\Hcal_\tau$, it follows that any simple object $S$ in $\Hcal_\tau$ is either of the form $S = F$ for some $F$ in $\Fcal$ or $S=T[-1]$ for some $T$ in $\Tcal$.  In other words, the simple objects in $\Hcal_\tau$ correspond to certain \emph{modules} in $\Mod{R}$.  The aim of this section is to identify these modules. We remark that our results remain valid when replacing $\Mod{R}$ by an arbitrary abelian category. 

\begin{definition}\label{def: tf/t t/tf}
Let $\tau = (\Tcal, \Fcal)$ be a torsion pair in $\Mod{R}$.

\noindent A non-zero module $T$ is called \textbf{almost torsionfree} if the following conditions are satisfied.\begin{enumerate}
\item[(ATF1)] Every proper submodule of $T$ is contained in $\Fcal$.
\item[(ATF2)] For each short exact sequence $0 \to A \to B \to T \to 0$, if $B$ is in $\Tcal$, then $A$ is in $\Tcal$.
\end{enumerate}
A non-zero module $F$ is called \textbf{almost torsion} if the following conditions are satisfied.\begin{enumerate}
\item[(AT1)] Every proper quotient of $F$ is contained in $\Tcal$.
\item[(AT2)] For each short exact sequence $0 \to F \to A \to B \to 0$, if $A$ is in $\Fcal$, then $B$ is in $\Fcal$.
\end{enumerate}
\end{definition}

\begin{remark} Any torsionfree module is trivially almost torsionfree and any torsion module is trivially almost torsion. The condition (ATF1) implies that if an almost torsionfree object is not torsionfree, then it must be torsion.  Similarly, any almost torsion object is either torsion or torsionfree. We will consider the non-trivial cases: the objects contained in $\Tcal$ that are almost torsionfree and the objects contained in $\Fcal$ that are almost torsion.  These objects are also known as \textbf{torsion, almost torsionfree} and \textbf{torsionfree, almost torsion} respectively. \end{remark}

\begin{example} \label{ex:ht pair}
Suppose that the torsion pair $\tau = (\Tcal, \Fcal)$ in $\Mod{R}$ is hereditary.   If $T$ is a torsion, almost torsionfree module, then (ATF1) implies that $T$ is simple. Conversely, if $T \in \Tcal$ is simple, then (ATF1) is clearly satisfied, and (ATF2) follows from
the hereditary property of $\Tcal.$

Next, we show that the torsionfree, almost torsion modules are precisely the modules in $\Cogen^2 (\Inj (\Fcal)) \isom \Mod{R}/\Tcal$ which become simple in the localisation.
To see that 
a torsionfree, almost torsion module $F$ belongs to  
$\Cogen^2 (\Inj (\Fcal))$, take the injective envelope of $F,$ 
\begin{equation} \label{Eq:tf/at envelope}
\xymatrix@1{0 \, \ar[r] & \, F \, \ar[r]^-{a} & \, E(F) \, \ar[r] & \, \Omega^{-1}(F) \, \ar[r] & \, 0,}
\end{equation} which is torsionfree by the hereditary property.
Condition (AT2) implies that $\Omega^{-1}(F)$ too is torsionfree; if we take its injective envelope, we get a copresentation of $F = F_{\Tcal}$ by torsionfree injective modules. 
Conversely, any module $F$ in $\Cogen^2 (\Inj (\Fcal))$ 
satisfies condition (AT2). For, suppose w.l.o.g. that there is a short exact sequence $0 \to F \to A \to B \to 0$ with  $A$  in $\Fcal$ and $B \neq 0$ in $\Tcal$. Then we  have a commutative diagram with exact rows \begin{equation} \xymatrix@1{0 \, \ar[r] & \, F \, \ar@{=}[d]\ar[r] & \, A \,\ar[d]\ar[r] & \, B\,\ar[d]^h \ar[r] & \, 0,\\ 0 \, \ar[r] & \, F \, \ar[r] & \, E(F) \, \ar[r] & \, \Omega^{-1}(F) \, \ar[r] & \, 0,}\end{equation} 
where $\Omega^{-1}(F)\in \Fcal$ and thus $h=0$. But then the upper row is split exact, a contradiction. 
Now it is easy to see that a module  in $\Cogen^2 (\Inj (\Fcal)),$ regarded as an object of the localisation, contains no proper subobjects, and must therefore be simple, if and only if it satisfies condition (AT1).

Finally, observe that the torsionfree, almost torsion module $F$ is uniform, for if $M_1 \cap M_2 = 0$ are two nonzero submodules of $F,$ then, by (AT1), the direct sum $F/M_1 \oplus F/M_2$ is a torsion module. As $F$ embeds in a canonical way into this direct sum, the hereditary property would give the contradiction that $F$ was also torsion. We conclude that $E(F)$ is indecomposable and, because $\Omega^{-1}(F)$ also belongs to $\Fcal,$ the injective envelope $a \colon F \to E(F)$ in the short exact sequence~(\ref{Eq:tf/at envelope}) is a
special $\Fcal^{\perp_1}$-envelope of $F$ in $\Mod{R}$ (cf.\ Theorem~\ref{Thm: las and specials}(2)).
\end{example}

The following proposition is essentially a rephrasing of \cite[Lem.~2.3]{sentieri:2020}.

\begin{proposition}\label{prop: alternative ttf, tft}
Let $\tau = (\Tcal, \Fcal)$ be a torsion pair in $\Mod{R}$.\begin{enumerate}
\item The following statements are equivalent for a non-zero module $T$.
\begin{enumerate}
\item T is almost torsionfree.
\item For every exact sequence $0 \to X \to Y \overset{g}{\to} T$ with $X$ in $\Fcal$, either \begin{enumerate} \item $Y$ is in $\Fcal$; or \item $g$ is a split epimorphism. \bigskip\end{enumerate}
\end{enumerate}

\item The following statements are equivalent for a non-zero module $F$.
\begin{enumerate}
\item F is almost torsion.
\item For every exact sequence $F \overset{g}{\to} X {\to} Y \to 0$ with $Y$ in $\Tcal$, either \begin{enumerate} \item $X$ is in $\Tcal$; or \item $g$ is a split monomorphism. \end{enumerate}
\end{enumerate}
\end{enumerate}
\end{proposition}
\begin{proof}
We will prove (1), the argument for (2) is completely dual.

\textbf{(1) [(a)$\Rightarrow$(b)]} Assume $T$ is almost torsionfree and consider an arbitrary exact sequence $0 \to X \to Y \overset{g}{\to} T$ with $X$ in $\Fcal$. The case where $g$ is an epimorphism is covered by the dual of \cite[Lem.~2.3]{sentieri:2020} (noting that the argument does not require $T$ to be in $\Tcal$).  It remains to consider the case where $g$ is not an epimorphism.  Then $\Im(g)$ is in $\Fcal$ by (ATF1), so $Y$ is in $\Fcal$ because $\Fcal$ is closed under extensions.  

\textbf{[(b)$\Rightarrow$(a)]}  Suppose $T$ satisfies (1)(b).  By the dual of \cite[Lem.~2.3]{sentieri:2020} (noting again that the argument does not require $T$ to be in $\Tcal$), it suffices to show that every proper subobject of $T$ is in $\Fcal$.  But this follows immediately if we consider the exact sequence $0 \to 0 \to Y \to T$.
\end{proof}

\begin{remark}\label{brick labels} Almost torsionfree and almost torsion modules are closely related to the minimal (co)extending modules over finite-dimensional algebras introduced in \cite{barnard:carroll:zhu:2019}, and also the brick labelling given in \cite{Asai:2020} for functorially finite torsion pairs and in \cite{demonet:iyama:reading:reiten:thomas:2017} for general torsion pairs.  The precise connections between these concepts are made clear in \cite{sentieri:2020}. 

If $\tau=(\Qcal, \Ccal)$ is a cotilting torsion pair, then $\Ccal$ is closed under direct limits,  hence all torsion, almost torsionfree modules are finitely generated. On the other hand,  there may be  torsionfree, almost torsion modules which are not finitely generated, as Example~\ref{large tf/t} will show. \end{remark}

\begin{theorem}\label{thm: simples t/tf tf/t}
Let $\tau = (\Tcal, \Fcal)$ be a torsion pair in $\Mod{R}$.  The simple objects $S$ in the heart $\Hcal_\tau$ of the HRS-tilt of $(\Tcal, \Fcal)$ are precisely those of the form $S=T[-1]$ with $T$ torsion, almost torsionfree and $S = F$ with $F$ torsionfree, almost torsion.
\end{theorem}
\begin{proof}
Using the canonical exact sequence $0\to F\to S\to T[-1]\to 0$ in $\Hcal_\tau$ with $F\in\Fcal$ and $T\in \Tcal$, we see that a simple object $S$ is either of the form $S=F$ or $S=T[-1]$.  Let us show that an object  of the form  $S=F$ with $F\in\Fcal$ is simple if and only if $F$ is almost torsion.  The other case is proven dually. 

For the only-if part, we start by considering  a proper submodule $U$ of $F$. Then, since $F=S$ is simple in $\Hcal_\tau$, the map $h:U\to F$ gives rise to an epimorphism $h:U\to F=S$ in $\Hcal_\tau$, hence the module $F/U=\Coker(h)$ is contained in $\Tcal$ by Lemma \ref{heartlemma}.  So, (AT1) is verified.
To prove (AT2),  we consider a short exact sequence $0\to F\stackrel{h}{\to} A\to B\to 0$ in $\Mod{R}$ with $A\in\Fcal$.  Here $h:F= S\to A$ is a monomorphism in $\Hcal_\tau$, and so the module $B=\Coker(h)$ is contained in $\Fcal$ again by Lemma \ref{heartlemma}.

Conversely, we show that (AT1) and (AT2) imply that $S = F$ is simple. To this end, we  claim that every morphism $0\neq f:S\to A$ in $\Hcal_\tau$ is a monomorphism. Since $f$ factors through the torsion part of $A$ with respect to the torsion pair $(\Fcal,\Tcal[-1])$, we can assume that $A=C$ for some $C\in\Fcal$. Then $f$ is a morphism in $\Mod{R}$, and $f$ is a monomorphism (by (AT1)) with cokernel in $\Fcal$ (by (AT2)). But then it follows from Lemma \ref{heartlemma} that $\Ker_{\Hcal_\tau}(f)=0$, and the claim is proven.
\end{proof}

\begin{corollary}\label{cor: bricks}
Let $\tau = (\Tcal, \Fcal)$ be a torsion pair in $\Mod{R}$.
\begin{enumerate}
\item Assume that $T,T'$ are both torsion, almost torsionfree. If $g \colon T \to T'$ is non-zero, then $g$ is an isomorphism. 

\item Assume that $F,F'$ are both torsionfree, almost torsion. If $f \colon F \to F'$ is non-zero, then $f$ is an isomorphism. 
\end{enumerate}
\end{corollary}

In particular, we have shown that the torsion, almost torsionfree modules and the torsionfree, almost torsion modules are \textbf{bricks} (i.e., their endomorphism rings are division rings).

\section{Injective envelopes in a cotilting heart}\label{sec: inj env in heart}

In this section we will consider the case where our torsion pair $\tau = (\Qcal, \Ccal)$ is a cotilting torsion pair.  We know from \cite{colpi:gregorio:mantese:2007} that the associated cotilting heart $\Hcal_\tau$ is a Grothendieck category and so, in particular, has enough injectives.  Next we relate the injective envelopes of simple objects in  $\Hcal_\tau$ to the special approximation sequences induced by the perfect cotorsion pair $(\Ccal, \Ccal^{\perp_1}) = ({}^{\perp_1}C, ({}^{\perp_1}C)^{\perp_1})$.

\begin{proposition}\label{prop: inj envelope}
Let $\tau = (\Qcal, \Ccal)$ be a cotilting torsion pair with associated cotilting module $C$ and cotilting heart $\Hcal_\tau$.\begin{enumerate}
\item Let $M \in \Qcal$ and consider a short exact sequence $0\to  X \overset{a}{\to} Y \overset{b}{\to} M \to 0$ in $\Mod{R}$.  Let \[X \overset{a}{\to} Y \overset{b}{\to}  M\overset{c}{\to} X[1]\] be the corresponding triangle in $\Der{R}$.  The following statements are equivalent. \begin{enumerate} \item The morphism $b \colon Y \to M$ is a special $\Ccal$-cover of $M$ in $\Mod{R}$. 
\item The morphism $c[-1] \colon M[-1]{\to} X$ is an injective envelope of  $M[-1]$ in $\Hcal_\tau$.
\end{enumerate}

\item Let $M \in \Ccal$ and consider a short exact sequence $0\to M\overset{a}{\to} X\overset{b}{\to} Y\to 0$ in $\Mod{R}$. The following statements are equivalent.  \begin{enumerate}
\item The morphism $a \colon M \to X$ is a special $\Ccal^{\perp_1}$-envelope of $M$ in $\Mod{R}$.
\item The morphism $a \colon M{\to} X$ is an injective envelope of  $M$ in $\Hcal_\tau$.
\end{enumerate}
\end{enumerate}
\end{proposition}
\begin{proof}
\textbf{(1) [(a)$\Rightarrow$(b)]} Since $\Ccal$ is closed under submodules and $b$ is a special $\Ccal$-cover, $X\in\Ccal\cap\Ccal^{\perp_1}=\Prod(C)$, so $X$ is injective in $\Hcal_\tau$.  Moreover, it follows from Lemma \ref{heartlemma} that  there is an exact sequence \[0\to M[-1]\stackrel{c[-1]}{\to} X\stackrel{a}{\to} Y\to 0\] in $\Hcal_\tau$. 

It remains to check that $c[-1]$ is left minimal.
Consider an endomorphism $h\in{\rm End}_{\Hcal_\tau}(X)$ with $h\circ c[-1]=c[-1]$. Then there is $g\in{\rm End}_{\Hcal_\tau} (Y)$ yielding a commutative diagram whose rows are given by triangles

\[ \xymatrix{ M[-1] \ar[r]^{c[-1]} \ar@{=}[d] & X \ar[d]^{h} \ar[r]^{a} & Y \ar[d]^{g} \ar[r]^{b} & M \ar@{=}[d] \\
M[-1] \ar[r]^{c[-1]} & X \ar[r]^{a} & Y \ar[r]^{b}  & M
}\]  It follows that $b = b \circ g$ and hence $g$ is an isomorphism by the minimality of $b$.  As $g$ is an isomorphism, we conclude that $h$ is an isomorphism as desired.

\textbf{[(b)$\Rightarrow$(a)]}  Let $X' \overset{a'}{\to} Y' \overset{b'}{\to}  M\overset{c'}{\to} X'[1]$ be the triangle in $\Der{R}$ induced by a special $\Ccal$-cover $b'$ of $M$.  We have already seen that $c'[-1]$ is an injective envelope of $M[-1]$ in $\Hcal_\tau$.  Since injective envelopes are unique up to isomorphism, there exists an isomorphism $h \colon X \to X'$ and a commutative diagram:

\[ \xymatrix{ X \ar[r]^a \ar[d]^h_\cong & Y \ar@{.>}[d]^{f} \ar[r]^{b} & M \ar@{=}[d] \ar[r]^{c} & X[1] \ar[d]^{h[1]}_\cong \\
X' \ar[r]^{a'} & Y' \ar[r]^{b'} & M \ar[r]^{c'}  & X'[1]
}\] where the induced morphism $f$ must also be an isomosrphism.  It follows that $b$ is a $\Ccal$-cover of $M$ in $\Mod{R}$.
\smallskip

\textbf{(2) [(a)$\Rightarrow$(b)]} Since $M$ and $Y$ are in $\Ccal$, we have  $X\in\Ccal \cap \Ccal^{\perp_1}=\Prod(C)$, so $X$ is injective in $\Hcal_\tau$.  Moreover, it follows from Lemma \ref{heartlemma} that there is an exact sequence  $0\to M\stackrel{a}{\to} X\to Y\to 0$ in $\Hcal_\tau$. Finally,  $a$ is left minimal in $\Hcal_\tau$ since so is $a$ in $\Mod{R}$.

\textbf{[(b)$\Rightarrow$(a)]}  Let $a' \colon M \to X'$ be a special $\Ccal^{\perp_1}$-envelope of $M$.  We have seen that $a'$ is an injective envelope of $M$ in $\Hcal_\tau$.  Since injective envelopes are unique up to isomorphism,  there exists an isomorphism $h \colon X \to X'$ such that $ha = a'$.  It follows that $a$ is a special $\Ccal^{\perp_1}$-envelope of $M$.
\end{proof}

\begin{theorem}\label{Thm: las and specials}
Let $\tau = (\Qcal, \Ccal)$ be a cotilting torsion pair with associated cotilting module $C$.  Consider a short exact sequence $0 \to L \overset{a}{\to} M \overset{b}{\to} N \to 0$ in $\Mod{R}$.\begin{enumerate}
\item The following statements are equivalent. \begin{enumerate}
\item The module $N$ is torsion, almost torsionfree and the morphism $b$ is a special $\Ccal$-cover of $N$ in $\Mod{R}$.
\item The module $L$ is in $\Prod(C)$ and the morphism $a$ is a strong left almost split morphism in $\Ccal$.
\end{enumerate}
\item The following statements are equivalent. \begin{enumerate}
\item The module $L$ is torsionfree, almost torsion and the morphism $a$ is a special $\Ccal^{\perp_1}$-envelope of $L$ in $\Mod{R}$.
\item The module $M$ is in $\Prod(C)$ and the morphism $b$ is a strong left almost split morphism in $\Ccal$.
\end{enumerate}
\end{enumerate}
\end{theorem}
\begin{proof}
\textbf{(1)[(a)$\Rightarrow$(b)]}  By Theorem \ref{thm: simples t/tf tf/t}, the object $N$ is simple in the heart $\Hcal_\tau$ and by Proposition \ref{prop: inj envelope}, we have that $c[-1] \colon N[-1]\to L$ is an injective envelope, where $L \to M \to N \overset{c}{\to} L[1]$ is the completion of the exact sequence to a triangle in $\Der{R}$.  In particular, this means that $L$ is injective in $\Hcal_\tau$ and hence $L \in \Prod(C)$.  It follows from Lemma \ref{heartlemma}, that $0\to N[-1] \overset{c[-1]}{\to} L \overset{a}{\to} M \to 0$ is a short exact sequence in $\Hcal_\tau$.  Then Remark \ref{rem: las in cotilting class} tells us that $a$ is a strong left almost split morphism in $\Ccal$.

\textbf{[(b)$\Rightarrow$(a)]}  It follows from our assumptions that $L$ is injective in $\Hcal_\tau$. By Proposition \ref{prop: las inj env simple} and Remark \ref{rem: las in cotilting class}, the kernel $S:= \Ker_{\Hcal_\tau}(a)$ is simple and the inclusion $0\to S \overset{c}{\to} L$ is an injective envelope.  Moreover, since strong left almost split morphisms starting at an object are unique up to isomorphism,  $a$ is an epimorphism in $\Hcal_\tau$.  By Lemma \ref{heartlemma}, we have that $N$ is in $\Qcal$ and also that $N[-1]\cong S$.  By Theorem \ref{thm: simples t/tf tf/t} we have that $N$ is almost torsionfree.  Finally, since $c$ is an injective envelope, it follows from Proposition \ref{prop: inj envelope} that $b$ is a $\Ccal$-cover of $N$.

\textbf{(2)[(a)$\Rightarrow$(b)]}  By Theorem \ref{thm: simples t/tf tf/t}, we have that $L$ is simple in $\Hcal_\tau$ and by Proposition \ref{prop: inj envelope} the morphism $a$ is an injective envelope in $\Hcal_\tau$.  In particular, we have that $M$ is contained in $\Prod(C)$ and the sequence $0 \to L \overset{a}{\to} M \overset{b}{\to} N \to 0$ is exact in $\Hcal_\tau$. By Remark \ref{rem: las in cotilting class} and Proposition \ref{prop: las inj env simple}(1) we have that $b$ is a strong left almost split morphism in $\Ccal$.  

\textbf{[(b)$\Rightarrow$(a)]}  Since $M$ is in $\Ccal$ and $\Ccal$ is closed under submodules, we have that $L$ is also in $\Ccal$.  Therefore $0 \to L \overset{a}{\to} M \overset{b}{\to} N \to 0$ is a short exact sequence in $\Hcal_\tau$.  By our assumption, we have that $b$ is a strong left almost split morphism in $\Ccal$ and so, by Remark \ref{rem: las in cotilting class} and Proposition \ref{prop: las inj env simple}(2), we have that $L$ is simple and $a$ is an injective envelope.  By Theorem \ref{thm: simples t/tf tf/t} and Proposition \ref{prop: inj envelope}, we have shown that condition (a) holds.
\end{proof}

Next we show that the strong left almost split morphisms arising in Theorem \ref{Thm: las and specials} are the only strong left almost split morphisms in $\Ccal$ with domain contained in $\Prod(C)$.  

\begin{lemma}\label{lem: strong las mono or epi}
Let $\Mcal$ be a full subcategory of $\Mod{R}$ that is closed under subobjects.  Then the following are equivalent for a module $M$ in $\Mcal$.\begin{enumerate}
\item There is a left almost split morphism $f \colon M \to \tilde{M}$ in $\Mcal$ that is not a monomorphism.
\item There is a left almost split morphism $g \colon M \to \bar{M}$ in $\Mcal$ that is an epimorphism.
\end{enumerate}
Moreover, a strong left almost split morphism is either a monomorphism or an epimorphism.
\end{lemma}
\begin{proof}
\textbf{[(2) $\Rightarrow$ (1)]} is trivially true.  To prove \textbf{[(1) $\Rightarrow$ (2)]}, observe that, if $f \colon M \to \tilde{M}$ is a left almost split morphism in $\Mcal$ that is not a monomorphism, then $g \colon M \to \Im(f)$ is a left almost split morphism in $\Mcal$ that is an epimorphism because $\Im(f) \in \Mcal$ by assumption.\smallskip

The final statements follows immediately because strong left almost split morphisms are unique up to isomorphism.
\end{proof}

\noindent Define \[\Ncal_C := \{ N \in \Prod(C) \mid \exists N \to \bar{N} \text{ a (strong) left almost split morphism in } \Ccal\}.\]  Note that, by Proposition \ref{prop: las inj env simple} and the subsequent corollaries, the set $\Ncal_C$ does not depend on whether we choose to include the word strong or not.  The previous lemma shows that $\Ncal_C$ is a disjoint union  $\Ncal_C=\Mcal_C \sqcup \Ecal_C$ where \[\Mcal_C := \{ L \in \Prod(C) \mid \exists L \to \bar{L} \text{ a strong left almost split monomorphism in } \Ccal \} \:\:\text{  and }\] \[\Ecal_C := \{M \in \Prod(C) \mid \exists M \to \bar{M} \text{ a strong left almost split epimorphism in } \Ccal \}.\]  Note that $\Mcal_C$ consists of the modules $L \in \Prod(C)$ arising in Theorem \ref{Thm: las and specials}(1) and $\Ecal_C$ consists of the modules $M \in \Prod(C)$ arising in Theorem \ref{Thm: las and specials}(2).

\begin{corollary}\label{cor: strong las = injective envelope of simples}
The following statements hold for a module $N$ in $\Prod(C)$.
\begin{enumerate} 
\item $N \in \mathcal{N}_C$ if and only if $N$ is the injective envelope of a simple object in $\Hcal_\tau$.  In this case $N$ is isomorphic to an indecomposable direct summand of any cotilting module that is equivalent to $C$.
\item $N\in\Mcal_C$ if and only if $N$ is the injective envelope of $T[-1]$ in $\Hcal_\tau$ where $T$ is a torsion, almost torsionfree module with respect to $\tau$.  In this case $T$ is the cokernel of the strong left almost split monomorphism $N\to\bar{N}$ in $\Ccal$.
\item $N\in\Ecal_C$ if and only if $N$ is the injective envelope of $F$ in $\Hcal_\tau$ where $F$ is a torsionfree, almost torsion module with respect to $\tau$.  In this case $F$ is the kernel of the strong left almost split epimorphism $N\to\bar{N}$ in $\Ccal$.

\end{enumerate}
\end{corollary}
\begin{proof}
(1) The first statement follows immediately from Corollary \ref{cor: las simples}. The latter statement follows from the fact that the injective envelope of a simple object in a Grothendieck category is indecomposable and the fact that the injective envelopes of simple objects arise as direct summands of any injective cogenerator, up to isomorphism. 

The statements (2) and (3) follow directly from Theorem \ref{Thm: las and specials} and Proposition \ref{prop: inj envelope}.
\end{proof}

In some cases, the heart $\Hcal_\tau$ turns out to be locally finitely generated and in this case we have the converse of the second part of Corollary \ref{cor: strong las = injective envelope of simples}.

\begin{corollary}\label{cor: min inj cog}
Suppose $\Hcal_\tau$ is locally finitely generated.  \begin{enumerate}
\item Let $D \in \Prod(C)$.  Every $N\in \Ncal_C$ is isomorphic to a direct summand of $D$ if and only if $D$ is a cotilting module that is equivalent to $C$.
\item 
 Let $\tilde{C}$ be a special $\Ccal^{\perp_1}$-envelope of $\bigoplus_{N\in \Ncal_C} N$.  Then $\tilde{C}$ is a cotilting module that is equivalent to $C$ and, moreover, $\tilde{C}$ is isomorphic to a direct summand of every other cotilting module that is equivalent to $C$.\end{enumerate}

\end{corollary}
\begin{proof}
We have already seen that every $N \in \Ncal_C$ arises as a direct summand of a cotilting module $D$ that is equivalent to $C$.  The rest of the corollary follows from the corresponding statements for locally finitely generated Grothendieck categories (see, for example, \cite[Prop.~3.17 and Cor.~3.18]{kussin:laking:2020}).  Note that the special $\Ccal^{\perp_1}$-envelope of $\bigoplus_{N\in \Ncal_C} N$ becomes the injective envelope of $\bigoplus_{N\in \Ncal_C} N$ in the heart by Proposition \ref{prop: inj envelope}.
\end{proof}

\begin{example}\label{locally fg}
There are some important cases where we know that the heart $\Hcal_\tau$ is locally finitely generated and so we may apply Corollary \ref{cor: min inj cog}. 
\begin{enumerate}
\item If $C$ is a cotilting module of cofinite type, then $\Hcal_\tau$ is the heart of a compactly generated t-structure by \cite[Lemma 3.7]{angeleri:hrbek:2017} and \cite[Thm.~2.3]{bravo:parra:2019}. It follows from \cite[Thm.~8.31]{saorin:stovicek:2020} that  $\Hcal_\tau$  is locally finitely presented and hence locally finitely generated.  
\item A cotilting module $C$ is an elementary cogenerator if and only if $\Hcal_\tau$ is locally coherent.  One implication follows from the description of the heart as a localisation of the functor category given in \cite{stovicek:2014}, the other is shown in \cite[Thm.~5.12]{laking:2018}.
\item \cite[Thm.~15.31]{goebel:trlifaj:2006}, \cite[Thm.~5.2]{saorin:2017} If $R$ is a right noetherian ring, then {\rm (1)} and {\rm (2)} apply, and the finitely presented objects in $\Hcal_\tau$  are precisely the objects which belong to the bounded derived category $\Dcal^b(\mod R)$.
\end{enumerate}
\end{example}

\section{Neg-isolated modules}\label{sec: inj env of simples}

Let $R$ be a ring and consider the category $(\Lmod{R}, \Ab)$ of additive functors from the category $\Lmod{R}$ of finitely presented left $R$-modules to the category $\Ab$ of abelian groups.  This functor category is a locally coherent Grothendieck category.    We will use the notation $(M,-) := \Hom_R (M,-)$ for the representable objects in $(\mod{R}, \Ab),$ $M \in \mod{R},$ and we will write $[F, G]$ to denote the set $\Hom_{(\Lmod{R}, \Ab)}(F, G)$ of natural transformations from $F$ to $G$.

\subsection{Finitely generated subfunctors of $\Ical$} \label{sec:fg subfunctors} 

In this section, we 
study the subfunctors of  the most important object of $(\mod{R}, \Ab),$ the forgetful functor $\Ical \colon \mod{R} \to \Ab.$ 
There is a
natural isomorphism $\xymatrix@1{\Ical  \;\ar[r]^-{\iota} & \; (R,-)}$ between $\Ical$ and the functor represented by $R_R.$ 
For $M \in \mod{R},$ the $M$-component is given by the morphism $\xymatrix@1{M \; \ar[r]^-{\iota_M} & \; (R,M),}$ $m \mapsto (f_m \colon 1 \mapsto m)$ in $\Ab$. 
This is actually a morphism of $R$-modules which induces an isomorphism of pointed $R$-modules
$\xymatrix@1{(M,m) \; \ar[r]^-{} & \; ((R,M), f_m)}$ for each  $m \in M$.
So if $\phi \subseteq \Ical$ is a subfunctor of the forgetful functor, then $m \in \phi (M)$ if and only if $f_m \in \phi (R,M).$ In this way we observe that the rule $\phi \mapsto \phi (R,-)$ is a bijective correspondence between the subfunctors of the forgetful functor and those of $(R,-).$

Associated to a pointed finitely presented module $(M,m),$ is the finitely generated subfunctor $\im (f_m,-) \subseteq (R,-)$ induced by 
the natural transformation $(f_m, -):\xymatrix@1{(M,-) \; \ar[r]^-{}  & \; (R,-).}$ 
It corresponds to the  subfunctor $H_{M,m}$ of the forgetful functor $\Ical$ which takes $N \in \mod{R}$ to the finite matrix subgroup (or pp-definable subgroup) $H_{M,m}(N)=\{h(m)\mid h\in\Hom_R(M,N)\}$.
On the other hand, if $\phi \subseteq (R,-)$ is a finitely generated subfunctor, then there is a natural transformation  $\eta \colon (M,-) \to (R,-)$ with image $\im \eta = \phi.$ By Yoneda's Lemma, there exists an $R$-linear morphism $f_m:\xymatrix@1{R \; \ar[r]^-{} & \; M}$ such that $\eta = (f_m,-).$ Thus, every finitely generated subfunctor of $\Ical$ arises from a pointed finitely presented module $(M,m)$ in this way.

\medskip

The finitely presented objects of $(\mod{R}, \Ab)$ admit a canonical extension to $\Mod{R}$ that respects direct limits. As all representable functors and their finitely generated subfunctors are finitely presented, this pertains to the finitely generated subfunctors $\phi \subseteq \Ical$ of the forgetful functor.

\begin{proposition} \label{prop:free real} Let $(M,m)$ be a pointed finitely presented module and $\phi=H_{M,m} \subseteq \Ical$. 
Then $(M,m)$ is a \textbf{free realisation} of $\phi,$ in the sense that $m\in\phi(M)$, and whenever $(N,n)$ is a pointed module with $n\in\phi(N)$,  then there exists a morphism $h \colon (M,m) \to (N,n)$ of pointed modules.
\end{proposition}

\begin{proof} The case when $(N,n)$ is finitely presented is clear by definition of $H_{M,m}$. For the general case, use the fact that $\phi$
respects direct limits, so there exists a pointed finitely presented module $(N',n')$ and a morphism $(N',n') \to (N,n)$ of pointed modules  with the property that $n'\in\phi(N')$. Now use again the definition of $\phi$.
\end{proof}

\subsection{The pp-type of a pointed module} 
 The finitely generated subobjects of $\Ical$ form a modular lattice. We will now use techniques from the model theory of modules~\cite{ziegler:1984,prest:1988,prest:2009} to investigate this lattice.
The model theoretic approach allows us to represent the finitely generated subfunctors of $\Ical$ by formulas $\phi (x)$ in a certain first-order language. The formula endows the corresponding subfunctor $\phi\subseteq \Ical$ with semantic content, so that if $(M,m)$ is a pointed right $R$-module, we may evaluate the statement $\phi (m)$ as true in $M,$ denoted by $M \models \phi (m),$ or not. In other words, the $M$-component of the inclusion $\phi \subseteq \Ical$ consists of the ``solutions'' in $M$ to the formula $\phi (x).$ We refer to the formulas $\phi (x)$ as \textbf{pp-formulas}, see \cite[\S 12.2]{prest:2009} for details.

The \textbf{pp-type,} denoted by $\pp (N,n),$ of a pointed module $(N,n)$ is the collection of pp-formulas $\phi (x),$ for which $N \models \phi (n).$ Equivalently, we can think of $\pp (N,n)$ as the collection of finitely generated subfunctors $\phi \subseteq \Ical$ for which $n \in \phi (N).$ As such, the pp-type $\pp (N,n)$ may be regarded as a filter $\Phi$ in the lattice of finitely generated subfunctors of the forgetful functor $\Ical \in (\mod{R}, \Ab);$ we say that $(N,n)$ is a \textbf{realisation} of $\Phi.$ The Completeness Theorem of first-order logic ensures that every filter $\Phi$ arises as the pp-type $\Phi = \pp (N,n)$ of some pointed module. The functorial property of pp-formulas ensures that if $f \colon (M,n) \to (N,n)$ is a morphism of pointed $R$-modules, then $\pp (M,m) \subseteq \pp (N,n).$

The general question thus arises of when an inclusion $\pp (M,m) \subseteq \pp (N,n)$ of pp-types is induced by a morphism of pointed modules. If $M$ is finitely presented, it is easy to see  that  $\pp (M,m)$ is the principal filter generated by the subfunctor $H_{M,m}\subseteq \Ical$. 
Proposition~\ref{prop:free real} therefore  implies that every inclusion $\pp (M,m) \subseteq \pp (N,n)$ is induced by a morphism $h \colon (M,m) \to (N,n).$ Next we consider a condition on the module $(N,n)$ that ensures the existence of a morphism of pointed modules.

\subsection{Pure-injective modules}\label{Sec: functor cats} 
The  assignment $M \mapsto (M \otimes_R -)$ defines a fully faithful right exact functor \begin{equation}\label{eq: tensor} \cY \colon \Mod{R} \to (\Lmod{R}, \Ab)\end{equation}   
which is called  coYoneda embedding. It allows us to consider the exact structure of the functor category inside the module category.  This is known as the pure exact structure in $\Mod{R}$.
 
\begin{definition}
A short exact sequence $\xymatrix@1@C=25pt{0 \, \ar[r] & \, N \, \ar[r]^f & \, M \, \ar[r]^g & \, L \, \ar[r] & \, 0}$ in $\Mod{R}$ is called \textbf{pure} if 
\begin{equation} \label{eq: pure exact} 
\xymatrix@1@C=25pt{0 \, \ar[r] & \, N \, \otimes_R - \, \ar[r]^-{f \tensor -} & \, M \otimes_R - \, \ar[r]^-{g \tensor -} & \, L \otimes_R - \, \ar[r] & \, 0}
\end{equation} 
is exact in $(\Lmod{R}, \Ab)$. In this case, we refer to $f$ as a \textbf{pure monomorphism} and to $g$ as a \textbf{pure epimorphism}.  
\end{definition}

It is natural to consider the modules that are injective with respect to the pure exact structure.

\begin{definition}
A module $N$ is called \textbf{pure-injective} if every pure exact sequence of the form (\ref{eq: pure exact}) is a split exact sequence.
\end{definition}

Clearly any module that becomes an injective object under the coYoneda embedding is pure-injective and, in fact, all injective objects in $(\Lmod{R}, \Ab)$ arise in this way. That is, the coYoneda embedding restricts to an equivalence 
$$\cY \colon \xymatrix@1{\Pinj(R) \, \ar[r]^-{\sim} & \, \Inj(\Lmod{R}, \Ab)}$$ where $\Pinj(R)$ denotes the full subcategory of pure-injective objects in $\Mod{R}$ and $\Inj(\Lmod{R}, \Ab)$ denotes the full subcategory of injective objects in $(\Lmod{R}, \Ab)$. Furthermore, if $M_R$ is an $R$-module with pure-injective envelope $\iota \colon M \to \PE (M),$ then the corresponding monomorphism $\iota \tensor - \colon M \tensor_R - \to \PE (M) \tensor_R -$ is the injective envelope of $M \tensor_R -$ in 
$(\Lmod{R}, \Ab).$ 

A pure monomorphism $f \colon M \to N$ may also be characterised~\cite[Proposition 2.1.6]{prest:2009} in terms of pp-formulas, which is equivalent to the condition that  $\pp (M,m) = \pp (N, f(m))$ for every $m \in M$. So if $\Phi$ is a filter of finitely generated subfunctors of the forgetful functor with realisation $\Phi = \pp (N,n),$ then the pure-injective envelope $\Phi = \pp (\PE (N), n)$ too is a realisation of $\Phi.$ 

\begin{remark} \label{rmk:pi hulls}~\cite[Corollary 3.3 (1)]{ziegler:1984} If $N$ is pure-injective, then every inclusion \linebreak $\pp (M,m) \subseteq \pp (N,n)$ is induced by a morphism $f \colon (M,m) \to (N,n).$ \end{remark}

\subsection{Definable subcategories of modules}\label{sec:pi adj}
We are interested in localisating the functor category $(\Lmod{R}, \Ab)$ at hereditary torsion classes associated to a particular kind of category of modules.

\begin{definition}
A full subcategory $\Dcal$ of $\Mod{R}$ is called \textbf{definable} if it is closed under products, pure submodules and directed colimits. 
\end{definition}

A definable subcategory $\Dcal \subseteq \Mod{R}$ is closed under pure-injective envelopes, so that its image under the coYoneda embedding $\cY (\Dcal) = \Dcal \tensor - \subseteq (\Lmod{R}, \Ab)$ is closed under injective envelopes in $(\Lmod{R}, \Ab).$ It follows that the torsion pair $(\Tcal_{\Dcal}, \Cogen (\Dcal \tensor -))$ in $(\Lmod{R}, \Ab)$ cogenerated by $\Dcal \tensor -,$ is hereditary. 

\begin{notation}\label{not: def loc}
We will denote the localisation $(\Lmod{R}, \Ab)/{\Tcal_\Dcal}$ by $(\Lmod{R}, \Ab)_\Dcal$ and the corresponding localisation functor by \[(-)_\Dcal \colon (\Lmod{R}, \Ab) \to (\Lmod{R}, \Ab)_\Dcal.\]  
We will denote the Hom-spaces between two objects $F,G$ in $(\Lmod{R}, \Ab)_\Dcal$ by $[F, G]_\Dcal$. 
\end{notation}

If $\Dcal \subseteq \Mod{R}$ is a definable subcategory, then the hereditary torsion pair $(\Tcal_{\Dcal}, \Cogen (\Dcal \tensor -))$ in $(\Lmod{R}, \Ab)$ is of finite type in the sense of the following definition. 
We refer the reader to \cite{herzog:1997} for more details on the theory surrounding hereditary torsion pairs of finite type in locally coherent Grothendieck categories. 

\begin{definition}
A torsion pair $(\Tcal, \Fcal)$ (not necessarily hereditary) in a Grothendieck category $\Gcal$ is said to be of \textbf{finite type} if the torsionfree class $\Fcal$ is closed under directed limits in $\Gcal.$
\end{definition} 

The following theorem will be very important in what follows. For details, we refer to \cite[\S 12.3]{prest:2009}.
\begin{theorem}\label{thm: def loc}
The rule $\Dcal \mapsto (\Tcal_{\Dcal}, \Fcal_{\Dcal})$ is a bijective correspondence between the collection of definable subcategories $\Dcal \subseteq \Mod{R}$ and hereditary torsion pairs in $(\Lmod{R}, \Ab)$ of finite type, with inverse given by 
$(\Tcal, \Fcal) \mapsto \Dcal = \{ M \in \Mod{R} \, | \, M \tensor_R - \in \Fcal \; \}.$
\end{theorem}

Consider the functor $ \Mod{R} \to (\Lmod{R}, \Ab)_\Dcal$ given by the composition of the functor (\ref{eq: tensor}) with the localisation functor $(-)_\Dcal$.  Since the injective objects in $(\Lmod{R}, \Ab)_\Dcal$ coincide with the injective objects in $(\Lmod{R}, \Ab)$ that are contained in the torsionfree class, this functor restricts to an equivalence of categories \[ \Pinj(\Dcal) \overset{\sim}{\rightarrow} \Inj\left((\Lmod{R}, \Ab)_\Dcal\right)\]  where $\Pinj(\Dcal)$ denotes the full subcategory of pure-injective objects in $\Dcal$ and $\Inj\left((\Lmod{R}, \Ab)_\Dcal\right)$ denotes the full subcategory of injective objects in $(\Lmod{R}, \Ab)_\Dcal$.

 \begin{remark}\label{rem: tf injectives} In the proofs below we will use the following observation several times. If $\mathcal{D}$ is a definable subcategory of $(\Lmod{R}, \Ab)$, then, for any $M \in \mathrm{Pinj}(\mathcal{D})$ and any $R$-module $L$, we have that \begin{equation}\label{hom} [\,(L\otimes -)_\mathcal{D}, (M\otimes -)_\mathcal{D}\,]_\mathcal{D} \cong [(L\otimes -), (M\otimes -)] \cong \mathrm{Hom}_R(L, M).\end{equation}  This follows directly from the fact that $(M\otimes -)$ is injective and torsionfree with respect to the torsion pair $(\Tcal_\Dcal, \Fcal_\Dcal)$ induced by $\mathcal{D}$ in $(\Lmod{R}, \Ab)$.  In particular, we have \begin{equation}\label{homs} [\,(R\otimes -)_\mathcal{D}, (M\otimes -)_\mathcal{D}\,]_\mathcal{D} \cong [(R\otimes -), (M\otimes -)] \cong \mathrm{Hom}_R(R, M) \cong M.\end{equation} \end{remark}

\subsection{Neg-isolated pure-injective modules}\label{subsec: inj env of simples}
Next we consider the pure-injective objects in a given definable subcategory that correspond to injective envelopes of simple objects in the Grothendieck category $(\Lmod{R}, \Ab)_\Dcal$. They will be characterised by a condition on pp-types.

A filter $\Phi$ in the lattice of finitely generated subfunctors of the forgetful functor  in $(\mod{R}, \Ab)$ will be called a \textbf{$\Dcal$-filter} if it admits a realisation $\pp (D,d) = \Phi$ with $D$ in $\Dcal.$ Recall that we can always choose $D$ to belong to the full subcategory $\Pinj(\Dcal)$ of pure-injective objects in $\Dcal$.

\begin{theorem}\label{prop: neg is env of simple}\label{prop: las in F = neg-isol}
Let $\Dcal$ be a definable subcategory in $\Mod{R}$.  The following statements are equivalent for an indecomposable pure-injective module $N$ in $\Dcal$.
\begin{enumerate}
\item $(N\otimes_R-)_\Dcal$ is the injective envelope of a simple object in $(\Lmod{R}, \Ab)_\Dcal$.
\item There exists a left almost split morphism $N \to N^+$ in  $\Pinj(\Dcal)$.
\item If $n\in N$ is a nonzero element and   $\Phi=\pp (N,n)$ is the associated pp-type, then there exists   a $\Dcal$-filter $\Phi^+ \supset \Phi$ which properly contains $\Phi$ such that whenever a $\Dcal$-filter $\Psi \supset \Phi$ properly contains $\Phi,$ then $\Psi \supseteq \Phi^+.$ 
\item $N$ is the source of an almost split morphism in $\Dcal$.
\end{enumerate}
\end{theorem}
\begin{proof}
Notice that, since $N$ is an indecomposable pure-injective module, the endomorphism ring of $N$ is then local. This implies that every endomorphism of pointed modules $ (N,n) \to (N,n)$ with $n\in N$ being nonzero is an automorphism. The equivalence \textbf{[(1) $\Leftrightarrow$ (2)]} follows directly from Corollary \ref{cor: inj env simp} and the discussion following Notation \ref{not: def loc}. 

\textbf{[(2) $\Rightarrow$ (3)]}
Let $n \in N$ be nonzero and let $f \colon N \to N^+$ be a left almost split morphism in $\Pinj (\Dcal).$ Set $n^+=f(n)$ and consider the $\Dcal$-filters $\Phi=\pp (N,n)$ and  $\Phi^+=\pp (N^+,n^+)$.
Clearly $\Phi^+\supseteq \Phi$, and equality would imply by Remark~\ref{rmk:pi hulls} that there is a map of pointed modules $(N^+,n^+)\to (N,n)$.
 Because $f$ is not a split monomorphism, we get that $\Phi^+\supset \Phi$ is strictly larger. 
On the other hand, if $\Psi \supset \Phi$ is any strictly larger $\Dcal$-filter with realisation $\Psi=(U,u)$ where $U$ in $\Pinj (\Dcal)$, then there is a morphism $g \colon (N,n) \to (U,u)$ which is not a split monomorphism. Then $g$ factors through $f,$ and $\pp (U, u) \supseteq \pp (N^+,n^+)$, that is, $\Psi \supseteq \Phi^+$.

\textbf{[(3) $\Rightarrow$ (4)]} Given a  nonzero element $n \in N$, set $\Phi=\pp (N,n)$ and pick $\psi\in\Phi^+\setminus\Phi$ with  a free realisation  $(M,m)$, i.e.~a finitely presented pointed module $(M,m)$ such that  $\psi=H_{M,m}$. Consider the pushout of $(N,n)$ with  $(M,m)$, postcomposed with a $\Dcal$-approximation (which exists e.g.~by~\cite[Proposition 3.4.39]{prest:2009}) as in the middle row of
$$\xymatrix@R=35pt@C=35pt{(R,1) \ar[r]^{f_m} \ar[d]^{f_n} \ar@{}[dr]|(.8){\ulcorner} & (M,m) \ar[d]^q \\
(N,n) \ar[r]^p \ar[d]^g & (K,k) \ar@{.>}[ld] \ar[r]^-{a_K} & (K_{\Dcal}, a_K (k)) =:(N^+,n^+)\ar@{.>}[lld]\\
(U,u)}$$
Clearly, $\pp (N^+, n^+)$ contains both $\Phi$ and $\psi$, so it must contain $\Phi^+$. Then the composition $a_K\circ p \colon (N,n) \to (N^+,n^+)$ cannot  be a split monomorphism. In fact, it must be a left almost split morphism in $\Dcal.$ For, suppose that we are given a morphism $g \colon N \to U$ in $\Dcal$ that is not a split monomorphism. Then there can't be a  morphism of pointed modules $(U,u)\to(N,n)$, and we infer from  Remark~\ref{rmk:pi hulls} that
$\pp (U,u)\supset \Phi$ is strictly larger. Thus $\psi \in \pp(U,u).$ By Proposition~\ref{prop:free real} there exist a morphism $(M,m) \to (U,u)$ from the free realisation of $\psi$ and a factorisation through the pushout. As $U \in \Dcal,$ this map from the pushout then factorises through its $\Dcal$-approximation, as required. 

\textbf{[(4) $\Rightarrow$ (2)]} Suppose there exists a left almost split morphism $h \colon N \to \bar{N}$ in $\Dcal$ and let $e \colon \bar{N} \to PE(\bar{N})$ be the pure-injective envelope of $\bar{N}$.  We will show that $g := eh$ is a left almost split morphism in $\Pinj(\Dcal)$.  Let $u \colon N \to U$ be a morphism in $\Pinj(\Dcal)$ that is not a split monomorphism.  Then there exists some $v \colon \bar{N} \to U$ such that $u = vh$.  Using that $e$ is a pure monomorphism and that $U$ is pure-injective, we have that there exists $f \colon PE(\bar{N}) \to U$ such that $fg = feh = vh = u$ as desired.
\end{proof}
The theorem above is  a relative version of results in~\cite[\S 5.3.5]{prest:2009} for the case $\Dcal=\Mod{R}$. Indeed, in that case condition (3)  means precisely that $\Phi=\pp (N,n)$ is a \textbf{neg-isolated pp-type}. This suggests the following terminology.
\begin{definition}
Let $\Dcal$ be a definable subcategory of $\Mod{R}$. A pure-injective module $N$ in $\Dcal$ is called \textbf{neg-isolated in $\Dcal$} if it satisfies the equivalent conditions of Theorem~\ref{prop: neg is env of simple}.
\end{definition}

We have seen that left almost split morphisms in a cotilting class are intimately related to the injective envelopes of simple objects in the cotilting heart.  Since cotilting classes $\Ccal$ are always definable subcategories, it is natural to ask how the results of Section \ref{sec: inj env in heart} are related to the neg-isolated modules in $\Ccal$.

\begin{corollary}\label{cor: inj env neg-iso}
Let $C$ be a cotilting module with  torsion pair $\tau = (\Qcal, \Ccal)$. The  $R$-modules that become injective envelopes of simple objects in $\Hcal_\tau$  are precisely the neg-isolated modules in $\Ccal$ which lie in  $\Prod(C)$.
\end{corollary}
\begin{proof}
This follows immediately from Theorem \ref{prop: las in F = neg-isol} and Corollary \ref{cor: strong las = injective envelope of simples}.
\end{proof}

%
%

\subsection{Critical modules}\label{sec: critical}

Let $(\Qcal, \Ccal)$ be a cotilting torsion pair in $\Mod{R}$.  The aim of this section is to investigate the neg-isolated modules in $\Ccal$ that are domains of left almost split morphisms in $\Ccal$ that are epimorphisms.  In Lemma \ref{lem: strong las mono or epi} we saw that these coincide with the domains of left almost split morphisms in $\Ccal$ that are not monomorphisms; these are the critical neg-isolated modules.

\begin{definition} \label{def:cr neg isolated}
Let $\Dcal$ be a definable subcategory of $\Mod{R}$.  We call a neg-isolated module $N$ in $\Dcal$ \textbf{critical} in $\Dcal$ if there exists an morphism $h \colon N \to N^+$ that is a left almost split morphism in $\Pinj(\Dcal)$ such that $h$ is not a monomorphism.
\end{definition}

When the definable subcategory is a torsionfree class, we have the following alternative characterisation of critical modules showing that they are exactly the neg-isolated modules such that the associated strong left almost split morphism is an epimorphism.

\begin{proposition}\label{prop: epi and critical}
Let $(\Tcal, \Fcal)$ be a torsion pair in $\Mod{R}$ such that $\Fcal$ is a definable subcategory.  The following statements are equivalent for a module $N$ in $\Fcal$.\begin{enumerate}
\item $N$ is a critical neg-isolated module in $\Fcal$.
\item There exists a left almost split morphism $f \colon N \to \bar{N}$ in $\Fcal$ that is an epimorphism.  In particular, $f$ is a strong left almost split morphism.
\end{enumerate}
\end{proposition}
\begin{proof}
\textbf{[(1) $\Rightarrow$ (2)]}  Let $h \colon N \to N^+$ be a left almost split morphism in $\Pinj(\Fcal)$ that is not a monomorphism.  Note that $\Im(h)$ is contained in $\Fcal$ because $\Fcal$ is closed under submodules.  We will show that $\bar{h} \colon N \to \Im(h)$ is a left almost split morphism in $\Fcal$.  Since $h$ is not a monomorphism, it follows that $\bar{h}$ is not a split monomorphism.  Suppose $u \colon N \to U$ is a morphism in $\Fcal$ that is not a split monomorphism and consider $eu \colon N \to PE(U)$ where $e$ is the pure-injective envelope of $U$. Then $eu$ can't be a split monomorphism.  Since $h$ is left almost split in $\Pinj(\Fcal)$ and $PE(U)$ lies in $\Fcal$ by \cite[Thm.~3.4.8]{prest:2009}, there exists a morphism $f \colon N^+ \to PE(U)$ such that $eu = fh$.  Let $k \colon \Ker(h) \to N$ be the kernel of $h$. Then $euk = fhk = 0$ and, moreover, $uk = 0$ because $e$ is a monomorphism.  Using that $\bar{h}$ is the cokernel of $k$, we conclude that there exists a unique morphism $g \colon \Im(h) \to U$ such that $u = g\bar{h}$.\smallskip

\textbf{[(2) $\Rightarrow$ (1)]}  Let $h \colon N \to \bar{N}$ be a left almost split morphism in $\Fcal$ that is an epimorphism.  Then $N$ must belong to $\Pinj(\Fcal)$, since otherwise the pure-injective envelope $e:N\to PE(N)$ would factor through $h$ and $h$ would be an isomorphism. In the proof of Theorem \ref{prop: las in F = neg-isol}, we showed that $eh \colon N \to \bar{N} \to PE(\bar{N})$ is a left almost split morphism in $\Pinj(\Fcal)$.  Since $h$ is an epimorphism, it is clear that $eh$ is not a monomorphism.
\end{proof}

We saw in Theorem~\ref{prop: neg is env of simple} that neg-isolated modules in definable subcategories are in bijection with injective envelopes of simple objects in the corresponding localisation of the functor category.  In the next lemma we identify which injective envelopes of simple objects give rise to critical modules.

\begin{lemma}
Let $\Dcal$ be a definable subcategory of $\Mod{R}$ and suppose $N$ is neg-isolated in $\Dcal$.  The following statements are equivalent. \begin{enumerate}
\item $N$ is critical in $\Dcal$.
\item $(N\otimes_R-)_\Dcal$ is the injective envelope of a simple object $S$ in $(\Lmod{R}, \Ab)_\Dcal$ such that $[(R\otimes_R-)_\Dcal, S]_\Dcal \neq 0$.
\end{enumerate}
\end{lemma}
\begin{proof}
Since $N$ is neg-isolated in $\Dcal$, we have the following set up according to Section \ref{Sec: functor cats}.  There exists a left almost split morphism $h \colon N \to N^+$ in $\Pinj(\Dcal)$ and, moreover, the morphism $(h\otimes -)_\Dcal \colon (N\otimes_R -)_\Dcal \to (N^+ \otimes -)_\Dcal$ is a left almost split morphism in $\Inj(\Lmod{R}, \Ab)_\Dcal$.  By Proposition \ref{prop: las inj env simple}, we see that the kernel of $(h\otimes-)_\Dcal$ is isomorphic to the monomorphism $i \colon S \to E(S) = (N\otimes_R-)_\Dcal$.

\textbf{[(1)$\Rightarrow$(2)]}  Suppose that $\Ker(h) \neq 0$.  Consider an element $k \in \Ker(h)$ and consider it as the morphism $k \colon R \to N$ that takes $1 \mapsto k$.  This yields a non-zero morphism $(k\otimes -) \in [(R\otimes_R-), (N\otimes_R -)]$, and it follows from Remark \ref{rem: tf injectives} that $(k\otimes -)_\Dcal \neq 0$.  Moreover, we have that $(h \otimes -)\circ(k\otimes-) = (h(k)\otimes-) = 0$ and so $(h \otimes -)_\Dcal\circ(k\otimes-)_\Dcal =0$.  Therefore, the morphism $(k\otimes-)_\Dcal$ factors through the kernel of $(h\otimes-)_\Dcal$. Thus we have a factorisation $i \circ g = (k\otimes -)_\Dcal$ and, in particular, we have a non-zero morphism $g \colon (R\otimes_R-)_\Dcal \to S$.

\textbf{[(2)$\Rightarrow$(1)]} Now suppose that there exists a non-zero morphism $f \colon (R\otimes_R-)_\Dcal \to S$ and consider $i \circ f \colon (R\otimes_R -)_\Dcal \to (N\otimes_R-)_\Dcal$.  Then, by Remark \ref{rem: tf injectives},  there exists a non-zero morphism $g \colon (R\otimes_R -) \to (N\otimes_R-)$ such that $g_\Dcal = i \circ f$.  Since the coYoneda embedding given in Section \ref{Sec: functor cats} is fully faithful, we can find a non-zero element $n \in N$ that determines a morphism $n \colon R \to N$ such that $1 \mapsto n$ and $(n \otimes -) = g$.  By definition, we have that $(n \otimes -)_\Dcal = i \circ f$ and so we have \[ (h(n) \otimes -)_\Dcal = ((h\otimes-)\circ(n\otimes-))_\Dcal = (h\otimes-)_\Dcal\circ(n\otimes-)_\Dcal = (h\otimes-)_\Dcal \circ i \circ f = 0.\]  By Remark \ref{rem: tf injectives}, this implies that $(h(n)\otimes-) = 0$ and hence $0 \neq n \in \Ker(h)$.
\end{proof}

This characterisation of critical modules in a definable subcategory $\Dcal$ allows us to show in the next proposition that they cogenerate $\Dcal$.  Moreover, it follows from this that the critical modules in a definable subcategory $\Dcal$ are related to the split injective modules in $\Dcal$.  A module $M$ in a full subcategory $\Mcal$ of $\Mod{R}$ is called \textbf{split injective in $\Mcal$} if every monomorphism $M \to M'$ with $M'$ in $\Mcal$ is a split monomorphism.

\begin{proposition}\label{prop: critical properties}
Let $\Dcal$ be a definable subcategory of $\Mod{R}$ and let $\Dcal_0$ denote the set of critical modules in $\Dcal$.\begin{enumerate}
\item The subcategory $\Dcal$ is contained in $\Cogen(\Dcal_0)$.
\item The split injective modules in $\Dcal$ are contained in the set $\Prod(\Dcal_0)$.
\item A module  $L$  belongs to $\Dcal_0$ if and only if it is neg-isolated and split injective in $\Dcal$.  
\end{enumerate}
\end{proposition}
\begin{proof}
\textbf{(1)} The definable subcategory $\mathcal{D}$ is closed under pure-injective envelopes and so it suffices to show that $ \mathrm{Pinj}(\mathcal{D}) \subseteq \Cogen(\mathcal{D}_0)$.  Let $N \in \mathrm{Pinj}(\mathcal{D})$ and let $n \in N$ and consider \[ 0 \to K \to (R\otimes -)_\mathcal{D} \overset{(n\otimes -)_\mathcal{D}}{\longrightarrow} (N\otimes -)_\mathcal{D}\]  Write $(n\otimes -)_\mathcal{D} = jq$ where $q \colon (R\otimes -)_\mathcal{D} \to (R\otimes -)_\mathcal{D}/K$ is the canonical quotient morphism and $j$ is a monomorphism.  Observe that the functor $(R\otimes -)\cong\Hom_R(R,-)$ is finitely presented, and so is its localisation $(R\otimes -)_\mathcal{D}$ by \cite[Prop.~2.15]{herzog:1997}. Then $K$ must be contained in a maximal subobject $M$, which induces an epimorphism $r \colon (R\otimes -)_\mathcal{D}/K \to (R\otimes -)_\mathcal{D}/M$ such that $p = rq$ where $p$ is the quotient morphism $(R\otimes- )_\mathcal{D} \to (R\otimes- )_\mathcal{D}/M$.  Since $S := (R\otimes  -)_\mathcal{D}/M$ is simple and there is a non-zero morphism $(R\otimes- )_\mathcal{D} \to S$, it follows that there is a critical module $L_n$ in $\mathcal{D}$ such that $i \colon S \to (L_n \otimes -)_\mathcal{D}$ is an injective envelope.
  
  As $(L_n \otimes -)_\mathcal{D}$ is injective, there exists a morphism $h \colon (N\otimes -)_\mathcal{D} \to (L_n\otimes -)_\mathcal{D}$ such that the following diagram is commutative.
  \[ \xymatrix{ (R\otimes- )_\mathcal{D} \ar[rd]_q \ar[rrd]^{(n\otimes -)_\mathcal{D}} \ar[rdd]_p & & \\
  & (R\otimes- )_\mathcal{D}/K \ar[r]_j \ar[d]^r & (N\otimes -)_\mathcal{D} \ar[d]^h \\
  & S \ar[r]_i & (L_n\otimes -)_\mathcal{D}
  }\]  The isomorphisms in Remark \ref{rem: tf injectives} yield that $h \cong (h_n\otimes -)_\mathcal{D}$ for some $h_n \colon N \to L_n$.  The commutativity of the above diagram yields that $0 \neq (h_n\otimes -)_\mathcal{D} \circ (n\otimes -)_\mathcal{D} \cong (h_n(n)\otimes -)_\mathcal{D}$ and it follows that $h_n(n) \neq 0$.
  
  Applying this argument to every element of $N$, yields a monomorphism $N \to \prod_{n\in N} L_n$ with components equal to $h_n$.\smallskip
  
  \textbf{(2)}  By (1), every split injective in $\mathcal{D}$ is contained in $\Cogen(\mathcal{D}_0)$ and hence in $\Prod(\mathcal{D}_0)$.  
  
  \textbf{(3)} Let $m: L \to N$ be a monomorphism in $\Mod{R}$ with $N\in\Dcal$ and $L\in \Dcal_0$.  If $m$ is not split then, the composition $em$ is not split, where $e \colon N \to \mathrm{PE}(N)$ is the pure-injective envelope of $N$.  But then $em$ must factor through the left almost split morphism $h \colon L\to L^+$ in $\Pinj(\Dcal)$.  But this is not possible because $h$ is not a monomorphism.
  
 Conversely, if $L$ is neg-isolated and split injective in $\Dcal$, then there exists a left almost split morphism $h:L\to L^+$ in $\Pinj(\Dcal)$, and $h$ cannot be a monomorphism.
\end{proof}

  \begin{proposition}\label{splitinj} Now suppose $(\mathcal{Q}, \mathcal{C})$ is a cotilting torsion pair and let $\Ccal_0$ denote the set of critical modules in $\Ccal$. We fix an injective cogenerator $I$ of $\Mod{R}$ with a special $\Ccal$-cover $$0\to C_1\to C_0\stackrel{g}{\to} I\to 0$$   Then $C_0$ is split injective in $\Ccal$ and $\Prod(C_0) = \Prod(\Ccal_0)$.
\end{proposition}
 \begin{proof} 
First we show that $C_0$ is split injective. Consider a monomorphism $h:C_0\to C\in\Ccal$. Since $C$ is cogenerated by $C_0$ and $I$ is injective, there are  a cardinal  $\alpha$ and maps $e:C\hookrightarrow C_0^\alpha$ and $f:C_0^\alpha{\to} I$ such that $feh=g$. As $g$ is a $\Ccal$-cover, there is also $f':C_0^{\alpha}{\to} C_0$ such that $f=gf'$. Now the right minimality of $g$ yields that $h$ is a split monomorphism.

Next we show that $\Prod(C_0) = \Prod(\Ccal_0)$.  By the first part of the proof and Proposition \ref{prop: critical properties}(2), we have that $C_0 \in \Prod(\Ccal_0)$ and hence $\Prod(C_0)\subseteq \Prod(\Ccal_0)$. As in \cite[Lem.~1.1]{angeleri:tonolo:trlifaj:2001}, we see that $C_0$ is a  cogenerator of $\Ccal$.  Therefore, by Proposition \ref{prop: critical properties}(3), we have that $\Ccal_0 \subseteq \Prod(C_0)$.  Thus we have that $\Prod(\Ccal_0) = \Prod(C_0)$.
\end{proof}

 The Proposition above shows in particular that every cotilting class contains critical neg-isolated modules and torsionfree, almost torsion modules. In contrast, we will see in Example~\ref{no t/tf}  that torsion, almost torsionfree modules need not exist. 
  Notice moreover that in general  the class of split injectives in $\Ccal$ is properly contained in $\Prod(\Ccal_0)$,  cf.~Example~\ref{not splitinj}.

\begin{corollary}\label{cor: crit inj env}
The set of critical neg-isolated modules in $\Ccal$ coincides with the set $\Ecal_C$ of modules $M \in \Prod(C)$ admitting a strong left almost split epimorphism $M \to \bar{M}$ in $\Ccal$.
\end{corollary}
\begin{proof}
The statement follows immediately from Proposition \ref{prop: epi and critical} and Proposition \ref{splitinj} since $\Prod(C_0) \subseteq \Prod(C)$.
\end{proof}

\subsection{Special modules}\label{sec: special}

Let $\tau = (\Qcal, \Ccal)$ be a cotilting torsion pair with cotilting module $C$. 
We saw in Corollary \ref{cor: strong las = injective envelope of simples} that the injective envelopes of simple objects in the heart $\Hcal_\tau$ are exactly the objects in the set \[\Ncal_C := \{ N \in \Prod(C) \mid \exists N \to \bar{N} \text{ a (strong) left almost split morphism in } \Ccal\}.\]  By Corollary \ref{cor: inj env neg-iso} the elements of $\Ncal_C$ are the neg-isolated modules in $\Ccal$ which belong to $\Prod(C)$.  Moreover, we showed in Lemma \ref{lem: strong las mono or epi} that $\Ncal_C = \Ecal_C \sqcup \Mcal_C$ where \[\Ecal_C := \{ M \in \Prod(C) \mid \exists M \to \bar{M} \text{ a strong left almost split epimorphism in } \Ccal \} \:\:\text{  and }\] \[\Mcal_C := \{L \in \Prod(C) \mid \exists L \to \bar{L} \text{ a strong left almost split monomorphism in } \Ccal \}.\] By Corollary \ref{cor: crit inj env} we have that the elements of $\Ecal_C$ are the critical neg-isolated modules in $\Ccal$.  In this section we will identify which of the non-critical neg-isolated modules in $\Ccal$ are contained in $\Mcal_C$.  In other words, we wish to determine the non-critical neg-isolated modules in $\Ccal$ which are contained in $\Prod(C)$.

\begin{proposition}\label{prop: special} Let $\tau = (\Qcal, \Ccal)$ be a cotilting torsion pair.  The following statements are equivalent for a module $N$ in $\Ccal$.
\begin{enumerate}
\item $N$ is contained in $\Prod(C)$ and is a non-critical neg-isolated module in $\Ccal$.
\item $N$ is contained in the set $\Mcal_C$ of modules in $\Prod(C)$ that admit a strong left almost split monomorphism $f \colon N\to \bar{N}$ in $\Ccal$.
\item There exists a strong left almost split monomorphism $f \colon N \to \bar{N}$ in $\Ccal$ such that the cokernel of $f$ is torsion, almost torsionfree.
\item There exists a left almost split morphism $g \colon N \to N^+$ in $\Pinj(\Ccal)$ such that the cokernel of $g$ is not contained in $\Ccal$.
\end{enumerate}
\end{proposition}
\begin{proof} 
\textbf{[(1)$\Rightarrow$(2)]}  By Theorem~\ref{prop: las in F = neg-isol} and Proposition~\ref{cor: crit inj env}, if $N\in\Prod(C)$ is non-critical neg-isolated, then $N$ is contained in $\Ncal_C \setminus \Ecal_C = \Mcal_C$. 
\newline\noindent\textbf{[(2)$\Rightarrow$(3)]} This follows from Theorem \ref{Thm: las and specials}(1).
\newline\noindent\textbf{[(3)$\Rightarrow$(4)]}  Let $f \colon N \to \bar{N}$ be as in (3) and let $T := \Coker f$. By the proof of Theorem~\ref{prop: las in F = neg-isol}, we have that the composition $g := ef$ is a left almost split morphism in $\Pinj(\Ccal)$ where $e \colon \bar{N} \to N^+ := \mathrm{PE}(\bar{N})$ is the pure-injective envelope of $\bar{N}$.  We show that $Z:=\Coker g$ is not contained in $\Ccal$.  Thus we have a commutative diagram: \[ \xymatrix{  0 \ar[r] & N  \ar@{=}[d] \ar[r]^f  & \bar{N} \ar[d]^e \ar[r]  & T  \ar[d]^k \ar[r]  & 0 \\ 0 \ar[r] & N \ar[r]^g & N^+  \ar[r] & Z  \ar[r] & 0 }\] where $h$ must be non-zero.  Indeed, if $h=0$, then $e$ factors through $g$ and $f$ is a split monomorphism, a contradiction.  Since $T\in\Qcal$, we conclude that $Z$ is not in $\Ccal = \Qcal^{\perp_0}$.
\newline\noindent\textbf{[(4)$\Rightarrow$(1)]}  Let $g \colon N \to N^+$ be as in (4). Note that, by definition, the module $N$ is neg-isolated in $\Ccal$. Set $Z:=\Coker g$ and consider the special $\Ccal$-cover $0 \to X \to Y \to Z\to 0$ of $Z$.  Note that $X\in \Ccal \cap \Ccal^{\perp_1} = \Prod(C) \subseteq \Pinj(\Ccal)$ because $X$ is a subobject of $Y\in \Ccal$.  Since $N^+$ is in $\Ccal$, we have the following commutative diagram: \[ \xymatrix{  0 \ar[r] & \Img g  \ar[d]^k \ar[r]^i  & N^+ \ar[d] \ar[r]^{}  & Z  \ar@{=}[d] \ar[r]  & 0 \\ 0 \ar[r] & X \ar[r] & Y  \ar[r] & Z  \ar[r] & 0 }\] where $g=ip$ is the canonical factorisation of $g$ through $\Img g$.  Then $kp$ must be a split monomorphism.  Indeed, if $kp$ is not a split monomorphism, then there exists a morphism $l \colon N^+ \to X$ such that $kp = lg = lip$.  Since $p$ is an epimorphism, it follows that $k=li$. A standard argument shows that the bottom sequence splits, which is not possible because $Z$ would then be isomorphic to a summand of $Y\in \Ccal$.  We have therefore shown that $N$ is contained in $\Prod(C)$ because it is isomorphic to a direct summand of $X\in \Prod(C)$.  Moreover, since $kp =lg$ is a monomorphism, so is $g$ and hence $N$ is not critical. 
\end{proof}

\begin{definition}
Let $\tau = (\Qcal, \Ccal)$ be a cotilting torsion pair.  A neg-isolated module $N$ in $\Ccal$ is called \textbf{special} if it satisfies the equivalent conditions of Proposition~\ref{prop: special}. \end{definition}


\begin{corollary}\label{C_0 and C_1}
 Let $I$ be  an injective cogenerator of $\Mod{R}$ with a special $\Ccal$-cover \begin{equation}\label{approx}0\to C_1\to C_0\stackrel{g}{\to} I\to 0.\end{equation}
\begin{enumerate}
\item Every critical neg-isolated module in $\Ccal$ is a direct summand of $C_0$.
\item Every special neg-isolated module in $\Ccal$ is a direct summand of $C_1$.
\end{enumerate}

\end{corollary}
\begin{proof}
To prove the corollary we will use the following general property of neg-isolated modules (see \cite[Prop.~9.29]{prest:1988}): if a neg-isolated module $N$ in a definable subcategory $\Dcal$ is a direct summand of $\prod_{i\in I} M_i$ where $M_i \in \Dcal$ for all $i\in I$, then $N$ is a direct summand of $M_i$ for some $i\in I$.  Since $C_0\oplus C_1$ is a cotilting module equivalent to $C$, it follows from Corollary \ref{cor: strong las = injective envelope of simples}(1) that every $N \in \Ncal_C$ is a direct summand of $C_0\oplus C_1$.  By Proposition \ref{prop: critical properties}(3), the neg-isolated summands of $C_0$ are exactly the critical ones and so the special neg-isolated modules are all direct summands of $C_1$ by \cite[Prop.~9.29]{prest:1988}.

\end{proof}

\begin{corollary} \label{C:minimal candidate}
If $\tau = (\Qcal, \Ccal)$ is a cotilting torsion pair with cotilting module $C,$ then $\PE (\oplus_{\Ncal_C} \, N)$ is isomorphic to a direct summand of $C.$
\end{corollary}

\begin{proof}
The argument used in the first part of the proof of Corollary~\ref{C_0 and C_1} implies that every $N \in \Ncal_C$ arises as a direct summand of $C.$ 
Equivalently, the functor $N \tensor_R -,$ regarded as an object in the localised functor category $(\Lmod{R}, \Ab)_\C,$
is a coproduct factor of the injective object $C \tensor_R -.$ Since $N$ is neg-isolated, the \textbf{socle} $\soc (N \tensor_R -)$ is a simple object in $(\Lmod{R}, \Ab)_\C.$ It follows that
$\coprod_{N \in \Ncal_C} \soc (N \tensor_R -) \subseteq \soc (C \tensor_R -)$
and therefore that the injective envelope $$E (\coprod_{N \in \Ncal_C} \soc (N \tensor_R -)) = E (\coprod_{\Ncal_C} \, (N \tensor_R -)) = E((\oplus_{\Ncal_C} \, N) \tensor_R -) = \PE (\oplus_{\Ncal_C} \, N) \tensor_R -$$ is a coproduct factor of $C \tensor_R -,$ as required.
\end{proof}

 \section{Examples}
 
In this section, we discuss some examples that illustrate our results. We also study the special case of cotilting modules induced by ring epimorphisms, where we establish some  interesting properties of special and critical neg-isolated modules.
 
 
  \subsection{The Kronecker algebra}\label{Kron}
    Let  $\Lambda$ be the Kronecker algebra, i.e.~the path algebra of the quiver
$ \bullet\rlap{\raise3pt\hbox{$\xrightarrow{}$}}
        \lower4pt\hbox{$\xrightarrow[]{}$}\,\bullet$ over an algebraically closed field $k$. It is well known that the category of finite dimensional indecomposable modules
 admits a canonical trisection $(\p,\tube,\q)$, where $\p$ and $\q$ denote the families given by the preprojective and preinjective modules, respectively, and $\tube=\bigcup_{x\in\mathbb X} \tube_x$ is the tubular family formed by the regular modules and indexed over the projective line  $\mathbb X=\mathbb P^1(k$). Given a simple regular module $S$, we denote by $S_\infty$ and $S_{-\infty}$ the Pr\"ufer and the adic module on the corresponding ray and  coray, respectively. Further, we denote by $G$ the generic module. Recall from \cite{ringel:1979} that $\End_RG$ is a  division ring, and $G$ is the unique (up to isomorphism) indecomposable module which has infinite length over $\Lambda$, but finite length  over its endomorphism ring.  Moreover, $G$ cogenerates the class $\Fcal=\tube^{\perp_0}$ of all torsionfree modules, it generates the class $\Dcal={}^{\perp_0}\tube$  of all divisible modules, and the intersection $\Fcal\cap\Dcal=\Add G$  is equivalent to the category of all modules over a simple artinian ring $Q$ which is obtained as universal localization of $\Lambda$ at $\tube$ and is Morita equivalent to   $\End_RG$. 

A complete description of the hearts of all cotilting $\Lambda$-modules is given in \cite{rapa:2019}. Let us focus on two special cases.

\begin{example}\label{large tf/t}
Consider the cotilting torsion pair $\tau=(\Tcal, \Fcal)$ generated by $\tube$. It is cogenerated by the cotilting module 
$$C=G\oplus \prod\{\text{all adic modules } S_{-\infty}\}.$$ 
The torsion, almost torsionfree modules coincide with the simple regular modules, while $G$ is the only torsionfree, almost torsion module. 
We  refer to \cite{rapa:2019} for the first statement and prove the second for the reader's convenience.

{\rm (AT1)} Let $g:G\to B$ be a proper epimorphism, and $0\to B'\to B\to \overline{B}\to 0$  the canonical exact sequence with $B'\in\Tcal$ and $\overline{B}\in\Fcal$. Then $\overline{B}$ lies in $\Fcal\cap\Dcal=\Add G$ , so $G\stackrel{g}{\to}B\to \overline{B}$ is a morphism in $\Add G$ and is therefore zero or a split monomorphism. It follows that $\overline{B}=0$ and  $B\in\Tcal$.

{\rm (AT2)} Let $0\to G\to B\to C\to 0$ be an exact sequence with $B\in\Fcal$. By applying the functor $\Hom_\Lambda(S,-)$ given by a simple regular module $S$, we obtain an exact sequence $$\Hom_\Lambda(S,B)\to \Hom_\Lambda(S,C)\to\Ext^1_\Lambda(S,G)\cong D\Hom_\Lambda(G,S)$$ where the first and third term are zero. Hence $C\in\Fcal$.

We have shown that $G$ is almost torsion.
For the uniqueness, observe that any other torsionfree, almost torsion module $X\in\Fcal$ is cogenerated by $G$, hence $\Hom_\Lambda(X,G)\not=0$, and $X\cong G$ by Corollary~\ref{cor: bricks}.

\smallskip

It follows that $G$ is simple injective in the heart $\Hcal_\tau$. Moreover, every simple regular module $S$ gives rise to a short exact sequence $$0\to S_{-\infty}\stackrel{a}{\lra}  S_{-\infty}\stackrel{b}{\lra} S\to 0$$ as in Theorem B, and to a minimal injective coresolution
$$0\to S[-1]{\lra} S_{-\infty}{\lra} S_{-\infty}\to 0$$   of the simple torsionfree object $S[-1]$ in $\Hcal_\tau$.
We infer that $G$ is the only critical neg-isolated module in $\Fcal$, and the special neg-isolated modules coincide with the adic modules. These are the neg-isolated modules in $\Fcal$ which belong to $\Prod(C)$. Observe that also the  modules in $\p$ are neg-isolated in $\Fcal$, see Theorem~\ref{prop: las in F = neg-isol}.

Finally, let us remark that $\Hcal_\tau$ is not hereditary. Indeed, it is shown in \cite[5.2]{stovicek:kerner:trlifaj:2011} that the heart of a torsion pair  is  hereditary only if the torsion pair splits. But if $S_x$ denotes the simple regular in the tube $\tube_x$, then by \cite[Prop. 5]{ringel:1998} there is a non-split exact sequence $$0\to \bigoplus_{x\in \mathbb X} S_x\to  \prod_{x\in \mathbb X} S_x\to G^{(\alpha)}\to 0$$ with $\bigoplus_{x\in \mathbb X} S_x\in\mathcal \Tcal$ and $G^{(\alpha)}\in\Fcal$. This shows that $(\Tcal,\Fcal)$ is not a split torsion pair.
\end{example}

We have just seen that   torsionfree, almost torsion modules may be infinite dimensional, while this is not possible for  torsion, almost torsionfree modules according to Remark~\ref{brick labels}. The next example, however, exhibits a cotilting module  without torsion, almost torsionfree modules. 
\begin{example}\label{no t/tf}
Consider now the torsion pair $\tau=(\mathcal Q,\C)$ in $\Mod \Lambda$ generated by the set $\q$. It is cogenerated by the cotilting module 
$$C=G\oplus \bigoplus\{\text{all Pr\"ufer modules } S_{\infty} \}.$$ 
The heart $\Hcal_\tau$ is  locally coherent and hereditary, and it is equivalent to the category   of quasi-coherent sheaves over $\mathbb X$. In particular, all simple objects in $\Hcal_\tau$  are torsion. In other words, there are no   
 torsion, almost torsionfree modules, and the torsionfree, almost torsion modules coincide with the simple regular modules. For details, we refer again to \cite{rapa:2019}.
 
Every simple regular module $S$ gives rise to a short exact sequence $$0\to S\stackrel{a}{\lra}  S_{\infty}\stackrel{b}{\lra}  S_{\infty} \to 0$$ as in Theorem A, which can also be regarded as the minimal injective coresolution
of the simple torsion object $S$ in $\Hcal_\tau$.
The critical neg-isolated modules in $\Ccal$ thus coincide with the Pr\"ufer modules, and there are no special neg-isolated modules.
\end{example}

The last example also shows that 
in general the class $\Prod{C_0}$ in Proposition \ref{prop: critical properties} does not coincide with the class of split injectives in $\Ccal$.  
\begin{example}\label{not splitinj} Let  $\tau=(\mathcal Q,\C)$ be as in Example~\ref{no t/tf}. Given an injective cogenerator $I$ of $\Mod\Lambda$, there is a short exact sequence \begin{equation}\label{approx2} 0\to C_1\to C_0\stackrel{g}{\to} I\to 0\end{equation} where $g$ is a $\Ccal$-cover,  $C_0$ is a direct sum of Pr\"ufer modules and $C_1$ is a direct sum of copies of $G$, see \cite[Thm.~7.1]{reiten:ringel:2006}.  Observe that $C_1$ is a direct summand of a product of copies of  $C_0$ by \cite[Prop.~3]{ringel:1998}, but it is not split injective in $\Ccal$ because the sequence (\ref{approx2}) is not split.
\end{example}

\subsection{Hereditary torsion pairs}
In this section, we assume that $\tau=(\mathcal Q,\mathcal C)$ is a \textbf{hereditary} cotilting torsion pair with cotilting module $C.$ Equivalently, the torsionfree class $\Ccal$ is closed under injective envelopes. If $V \in \Ccal$ is split injective, then its injective envelope $V \to E(V)$ is a monomorphism in $\Ccal,$ and must therefore be an isomorphism. We conclude that the subcategory of split injective objects of $\Ccal$ is given by the category $\Inj (\Mod{R}) \, \cap \, \Ccal$ of torsionfree injective modules, and that if $\Ccal_0 = \Ecal_C$ denotes the set of critical modules in $\Ccal$ as in Proposition~\ref{prop: critical properties}, then $$\Prod (\Ccal_0) = \Ccal \, \cap \, \Inj (\Mod{R}).$$

Given a module $M$, we denote by $\C(M)$  a $\C$-cover of $M$ and by $\Kcal(M)$ its kernel,
$$\xymatrix@1{0 \, \ar[r] & \, \Kcal (M) \, \ar[r] & \, \C(M) \, \ar[r] & \, M \, \ar[r] & 0.}$$ 
If $C$ is a cotilting module cogenerating $\C$, then $\Kcal(M)$ is an object of $\Prod (C)$ which is uniquely determined by $M$, up to isomorphism. If $F \in \Ccal$ is a simple torsionfree module, then its injective envelope $E(F)$ is also torsionfree and so serves as its own $\Ccal$-cover, but if $Q \in \Qcal$ is a {\em torsion} simple module, then its injective envelope $E(Q)$ is not torsionfree so its $\Ccal$-cover is given by the epimorphism in
\begin{equation} \label{eq:inj C-cover}
\xymatrix@1{0 \, \ar[r] & \, \Kcal (E(Q)) \, \ar[r] & \, \C(E(Q))\, \ar[r]^-{c} & \, E(Q) \, \ar[r] & 0,}
\end{equation} 
where $\Kcal (E(Q)) \neq 0.$

In Example~\ref{ex:ht pair}, we elaborated on the equivalence of Theorem~\ref{Thm: las and specials}(2) for a hereditary torsion pair, showing that  the torsionfree, almost torsion modules correspond to the simple objects of the localisation $\Mod{R}/\Qcal,$
and that their injective envelopes are the  critical neg-isolated indecomposable pure-injectives. The following elaborates on the equivalence given by Theorem~\ref{Thm: las and specials}(1).

\begin{theorem} \label{T:simple C-cover}
Suppose that $(\Qcal, \Ccal)$ is a hereditary cotilting torsion pair and $Q_R \in \Qcal$ is a torsion simple module. The $\Ccal$-cover of $Q$ is given by the pullback of (\ref{eq:inj C-cover}) along its injective envelope
$$\xymatrix{0 \, \ar[r] & \, \Kcal (E(Q)) \, \ar[r]^-{a} \ar@{=}[d] & \, F \, \ar[r]^-{b} \ar@{>->}[d] & \, Q \, \ar[r] \ar@{>->}[d] & 0  \\
0 \, \ar[r] & \, \Kcal (E(Q)) \, \ar[r] & \, \C(E(Q)) \, \ar[r]^-{c} & \, E(Q) \, \ar[r] & 0.
}$$
Consequently, $F = \C(Q)$, $\Kcal (Q) = \Kcal (E(Q))$ and $\C(E(Q)) = E(\C(Q)).$ Furthermore, $a \colon \Kcal (Q) \to \C(Q)$ is a left almost morphism in $\C$ and $\Kcal (Q)$ is a special neg-isolated indecomposable pure-injective. 
\end{theorem}

\begin{proof}
 As in the argument used in the proof of Proposition~\ref{splitinj}, the module $\C(E(Q))$ must be split injective and therefore, by the hereditary property, injective. Moreover, any indecomposable summand of $\C(E(Q))$ which does not intersect $\Kcal (E(Q))$ must be isomorphic to $E(Q)$, which is impossible as $Q\in \Qcal$.
We conclude that the kernel morphism of $c$ is an injective envelope.
 This also implies that the embedding of $F$ into $\C(E(Q))$ is an injective envelope.

Furthermore, it follows  that $\Kcal (E(Q))$ is indecomposable. For, suppose that $\Kcal (E(Q)) = K_1 \oplus K_2$ were a proper decomposition. Neither of the summands can be injective, since they are contained in the kernel of a $\C$-cover. The injective envelope of 
$\Kcal (E(Q))$ as well as its cosyzygy would then be decomposable, contradicting the fact that $E(Q)$ is not. 

Because $F$ is torsionfree and $\Kcal (E(Q))$ belongs to $\C^{\perp_1},$ $b \colon F \to Q$ is a special $\C$-precover. As such, it contains the $\C$-cover $\C(Q) \to Q$ as a direct summand, whose kernel would be a direct summand of 
$\Kcal (E(Q)).$ But $\Kcal (E(Q))$ is indecomposable, and $Q$ is not torsionfree, so the $\C$-cover of $Q$ must contain $\Kcal (E(Q)),$ and properly so. As $Q$ is simple, so we see that $b \colon F = \C(Q) \to Q$ is the $\C$-cover of $Q$ and 
$\Kcal (Q) = \Kcal (E(Q)).$ 

The last statement follows from Theorem~\ref{Thm: las and specials}(1).
\end{proof}

Theorem~\ref{T:simple C-cover} allows us to infer that  the module $\PE (\oplus_{\Ncal_C} \, N)$ of Corollary~\ref{C:minimal candidate} that arises as a summand of any cotilting module $C$ for $\C$ is itself cotilting. Because it is a summand of a cotilting module, it certainly satisfies the first two conditions of being one. To verify the third, consider the injective cogenerator $I = \prod_{\Sim (\Mod{R})} \, E(S)$ of $\Mod{R},$ where the index set runs over the set of all simple modules. It may be decomposed as
$$I = \prod_{Q \in \Sim (\Qcal)} \, E(Q) \;\;\; \oplus \prod_{F \in \Sim (\Ccal)} \, E(F),$$ where the index set has been partitioned into the simple torsion and simple torsionfree modules, respectively. 
Take the special $\Ccal$-precover of $I$ given by the product of the respective $\Ccal$-covers,
$$\xymatrix@1{0 \, \ar[r] & \, \prod_{Q \in \Sim(\Qcal)} \, \Kcal (Q) \, \ar[r] & \, \prod_{Q \in \Sim(\Qcal)} \, \C(E(Q)) \; \oplus \, \prod_{F \in \Sim(\Ccal)} \, E(F) \, \ar[r] & \, I \, \ar[r] & 0.}$$
Because every $\Kcal (Q)$ is neg-isolated, the kernel belongs to $\Prod (\PE (\oplus_{\Ncal_C} \, N)).$ On the other hand, the middle term is injective and therefore belongs to $\C \cap \Inj (\Mod{R}) = \Prod (\C_0).$


\begin{corollary} \label{minimal version}
If $(\Qcal, \Ccal)$ is a hereditary cotilting torsion pair with cotilting module $C,$ then the module $\PE (\oplus_{\Ncal_C} \, N) = \PE ((\oplus_{\Ecal_C} \, M) \; \oplus \; (\oplus_{\Mcal_C} \, L)),$ given more explicitly by 
$$\PE (\bigoplus_{Q \in \Sim (\Qcal)} \, \Kcal (Q) \;\;\; \oplus \bigoplus_{F \in \Sim (\Mod{R}/\Qcal)} \, E(F)),$$ 
is a $1$-cotilting module for $\Ccal$ that is isomorphic to a direct summand of $C.$
\end{corollary}

\subsection{Commutative noetherian rings}
  Let now $R$ be a commutative noetherian ring. It is shown in \cite{angeleri:pospisil:stovicek:trlifaj:2014} that the cotilting torsion pairs are precisely the hereditary  torsion pairs in $\M$ with $R$ being torsionfree. They are parametrized by the subsets $P\subset\Spec R$ that are closed under specialization and satisfy $\Ass R\cap P=\emptyset$. 
More precisely,   $P$ corresponds to the  hereditary torsion pair $(\mathcal T,\mathcal F)$ given by 
$$\mathcal T=\{M\in \M\mid \Supp{M}\subset P\}$$
$$\mathcal F=\{M\in\M\mid \Ass{M}\cap P=\emptyset\}$$
where $\mathcal T$ contains all $E(R/\mathfrak p)$ with $\mathfrak p\in P$, and $\mathcal F$ contains all $E(R/\mathfrak q)$ with $\mathfrak q\in Q=\Spec R\setminus P$.  

We denote by $\Max P$ the set of all  maximal ideals of $R$ which lie in $P$ and by $\Max Q=\{\mathfrak q\in Q\mid  V(\mathfrak q)\setminus\{\mathfrak q\}\subset P\}$ the set of all prime ideals which are maximal in  $Q$.
\begin{proposition}\label{commnoeth} The  module  
$$PE(\bigoplus_{\mathfrak m \in \Max P} \Kcal(R/\mathfrak m))\, \oplus \,\bigoplus_{\mathfrak q \in \Max Q} E(R/\mathfrak q)$$
is a cotilting module cogenerating  $\Fcal$ which is isomorphic to a direct summand of  any other cotilting module  cogenerating $\Fcal$.
\end{proposition}
\begin{proof}
The torsion simple modules are precisely the modules of the form  $R/\mathfrak m$ with $ \mathfrak m \in \Max P$.
Further, since the ring is noetherian, the direct sum $\bigoplus_{F \in \Sim (\Mod{R}/\Tcal)} \, E(F)$ is (pure-)injective. 
By Corollary \ref{minimal version}, it remains to show that  the modules $E(F)$ with $ F \in \Sim (\Mod{R}/\Tcal)$ are precisely the indecomposable injectives  $E(R/\mathfrak q)$ with $\mathfrak q \in \Max Q$, up to isomorphism.

In order to verify this, recall first that any such $E(R/\mathfrak q)$ is torsionfree and therefore  cogenerated by critical modules. In particular $\Hom_R(E(R/\mathfrak q),E(F))\not=0$ for some   $F\in \Sim (\Mod{R}/\Tcal)$. But  $E(F)$ is indecomposable injective and thus of the form $E(R/\mathfrak p)$ for some $\mathfrak p\in V(\mathfrak q)$. Moreover $\mathfrak p\in Q$ as $E(F)\in\Fcal$. It follows by assumption that  $\mathfrak p=\mathfrak q$, hence $E(R/\mathfrak q)\cong E(F)$.

For the converse implication, 
we employ \cite[Theorem 5.2]{stovicek:trlifaj:herbera:2014} which states that 
$$C=\Kcal(\bigoplus_{\mathfrak m \in \Max P} E(R/\mathfrak m))\, \oplus \,\bigoplus_{\mathfrak q \in \Max Q} E(R/\mathfrak q)$$
is a cotilting module with the stated properties.  It follows that  every $E(F)$ is isomorphic to a direct summand of $C$. Now observe that $\Kcal(\bigoplus_{\mathfrak m \in \Max P} E(R/\mathfrak m))$ cannot have injective summands by minimality, so $E(F)$ must be isomorphic to some $E(R/\mathfrak q)$ with $\mathfrak q \in \Max Q$. 
\end{proof}

We remark that the  modules in $\Sim (\Mod{R}/\Tcal)$ need not have the form $R/\mathfrak q$, as we are going to see next.

\begin{example}
Let $P=\Max {\mathbb Z}$ be the (specialization closed) set of all maximal ideals of $\mathbb Z$. Then $Q=\Max Q=\{0\}$, and $(\Tcal, \Fcal)$ is the  torsion pair formed by the torsion and torsionfree abelian groups, respectively. So
the only simple object in $\Mod{\mathbb Z}/\Tcal$ is $\mathbb Q$, which coincides with its  injective envelope and is also the only  module which is critical neg-isolated in $\Fcal$. 
\end{example}

\subsection{Minimal cotilting modules}
Throughout this section, we assume that $C$ is a cotilting $R$-module with torsion pair $\tau=(\Qcal, \Ccal)$ and that there is  an injective cogenerator $I$ of $\Mod{R}$ with a special $\Ccal$-cover \begin{equation}\label{approx3}0\to C_1\to C_0\stackrel{g}{\to} I\to 0\end{equation}
such that
\begin{enumerate}
\item[(M1)] The left perpendicular category $\Ycal={}^{\perp_{0,1}} C_1$ is closed under direct products;
\item[(M2)] $\Hom_R(C_0,C_1)=0$.
\end{enumerate}
Cotilting modules with this property are called \textbf{minimal} and we call the associated class $\Ccal$ a \textbf{minimal cotilting class}.

Observe that condition (M1) amounts to the fact that the inclusion functor $\Ycal\to \Mod R$ has a left adjoint and a right adjoint. This is equivalent to the existence of a ring epimorphism $\lambda:R\to S$ which is {\bf pseudoflat}, i.e.~$\Tor_1^R(S,S)=0$,  such that $\Ycal$ is the essential image of the functor $\lambda_*:\Mod S\to \Mod R$ given by restriction of scalars. 

In fact, minimal cotilting modules are closely related to pseudoflat  ring epimorphisms. Let $k$ be a commutative ring such that $R$ is a $k$-algebra,  and let  $D=\Hom_k(-,E)$ be the duality induced by an injective cogenerator of $\Mod k$. Any injective ring epimorphism  $\lambda:R\to S$ with the property that $\Cogen\, D(_RS)\subset {}^{\perp_1}D(_RS)$ is pseudoflat and induces a minimal cotilting right $R$-module $C=D(_RS)\oplus D(_RS/R)$, for which the sequence  (\ref{approx3}) can be chosen as
\begin{equation}\label{approx4}0\to  D(_RS/R)\to D(_RS)\stackrel{g}{\to} D(_RR)\to 0\end{equation}
We refer to \cite[Section 4]{angeleri:hrbek:2019} for details. In particular, we have 
\begin{theorem} \cite[Theorem 4.16, Corollaries 4.18 and 4.19]{angeleri:hrbek:2019} 
The map assigning to a ring epimorphism $\lambda:R\to S$  the  class $\Cogen\, D(_RS)$ yields a bijection between 
\begin{enumerate}
\item[(i)] epiclasses of injective  ring epimorphisms $\lambda:R\to S$ with $\Cogen\, D(_RS)\subset {}^{\perp_1}D(_RS)$,
\item[(ii)] minimal cotilting classes.
\end{enumerate}
Moreover, the set in  (i) equals the set of 
\begin{enumerate}
\item[(i')] epiclasses of injective pseudoflat ring epimorphisms,
\end{enumerate}
provided that $R$ has weak global dimension at most one or is a commutative noetherian ring.
\end{theorem}

\begin{example} {\rm\bf (1)} Let  $\Lambda$ be the Kronecker algebra over a field $k$. Consider the pseudoflat ring epimorphism $\lambda:\Lambda\to \Lambda_\tube$ given by universal localization at $\tube$. We have a short exact sequence
$$0\to \Lambda{\to} \Lambda_\tube\cong G\oplus G \to \bigoplus \{ \text{all Pr\"ufer modules } S_{\infty}\}\to 0$$
 in $\Lambda$-{\rm Mod}. Applying the duality $D=\Hom_k(-,k)$ we obtain a short exact sequence in $\Mod \Lambda$
$$0\to \prod \{ \text{all adic modules } S_{-\infty}\}\to G\oplus G\to D(\Lambda)\to 0$$
and a minimal cotilting module equivalent to the cotilting module  in Example~\ref{large tf/t}.

\smallskip

{\rm\bf (2)}  \cite[Example 4.10]{angeleri:hrbek:2019} The cotilting module in Example~\ref{no t/tf} is not minimal. In fact, it is the only non-minimal cotilting $\Lambda$-module up to equivalence.

\smallskip

{\rm\bf (3)} \cite[Example 4.22 (3)]{angeleri:hrbek:2019} Over a commutative noetherian ring of Krull dimension at most one,  all cotilting modules are minimal. 
\end{example}

Let us collect some properties of minimal cotilting modules.  In what follows we will make use of the notation introduced before Corollary \ref{cor: strong las = injective envelope of simples}, noting that the set of neg-isolated modules in $\Ccal$ that are contained in $\Prod(C)$ coincides with $\Ncal_C$, the set of critical neg-isolated modules coincides with $\Ecal_C$, and the set of special neg-isolated modules coincides with $\Mcal_C$.  See Section \ref{sec: special} for more details.

\begin{lemma}\label{no common}
Let  assumptions and notation be as above, and suppose that the heart $\Hcal_\tau$ is locally finitely generated.
\begin{enumerate}
\item $\Prod (C_0)=\Prod(\Ecal_C)$ and $\Prod (C_1)=\Prod(\Mcal_C)$.
\item An indecomposable module $X$ lies in $\Prod (C)$ if and only if it lies in $\Prod (C_0)$ or $\Prod (C_1)$, and not in both.
\item If $\Ccal$ is a proper subcategory of $\M$, then it contains special neg-isolated modules, and $\Qcal$ contains torsion, almost torsionfree modules.\end{enumerate}
\end{lemma}
\begin{proof}
  The first statement in (1)  and the inclusion $\Prod (C_1)\supseteq\Prod(\Mcal_C)$ are true in general, see Proposition~\ref{splitinj} and Corollary~\ref{C_0 and C_1}(2).
 
 When $\Hcal_\tau$ is locally finitely generated, then the direct product of all injective envelopes of simple objects is an injective cogenerator of  $\Hcal_\tau$. Hence $\Prod  (C)=\Prod (\Ncal_C)$ where $\Ncal_C=\Ecal_C\cup\Mcal_C$ is the set of all neg-isolated modules in $\Ccal$ which belong to  $\Prod (C)$.  So, every module $X\in \Prod (C)$ admits a split monomorphism   $\iota:X\hookrightarrow M_0\oplus M_1$ where $M_i$ is a direct product of modules in $\Ccal_i$ for $i=0,1$. Then we can write ${\rm id}_X= \pi_0\iota_0+\pi_1\iota_1$ where $\iota_i$ and $\pi_i$ are the components of $\iota$ and of its left inverse $\pi:M_0\oplus M_1\to X$. 
 
 Now, if $X\in\Prod (C_1)$, then it follows from condition (M2) that $C_0\in {}^{\perp_0}X$, hence by condition (M1) we have $\Prod(C_0)\subseteq{}^{\perp_0}X$, which implies that $\pi_0:M_0\to X$ vanishes. Thus ${\rm id}_X= \pi_1\iota_1$ and $X$ lies in $\Prod (\Mcal_C)$. This concludes the proof of statement (1).
 
 Moreover, if $X$ is an indecomposable module in $\Prod(C)$, then it is pure-injective and therefore has a local endomorphism ring. This shows that $\pi_i\iota_i$ must be an isomorphism for $i=0$ or $i=1$, that is, $X$ lies  in $\Prod (C_0)$ or $\Prod (C_1)$.
 Finally, $X$ can't lie  in both, because  $\Prod(C_0)\subset {}^{\perp_0}C_1$ by conditions (M1) and (M2). This proves statement (2).
 
 For statement (3), observe that the $\C$-cover  in (\ref{approx3}) can't be an isomorphism, so $C_1\not=0$ and $\C_1\not=\emptyset$. The existence of torsion, almost torsionfree modules then follows from Theorem B.
\end{proof}

\begin{proposition}\cite{sentieri:2022}\label{colimit of wide}
Let  assumptions and notation be as above, and suppose that $R$ is right artinian. Assume further that $\C$ is a proper subcategory of $\M$, and let $\Wcal$ be the class of all modules with a finite filtration by torsion, almost torsionfree modules. Then the class of all direct limits of modules in $\Wcal$ coincides with the left perpendicular category ${}^{\perp_{0,1}} C_0$  of the module $C_0$ in (\ref{approx3}).
\end{proposition}

\begin{theorem}\label{properties of specials and criticals}
Let $R$ be a right artinian ring and let $C$ be a minimal cotilting module with associated ring epimorphism $\lambda:R\to S$.
\begin{enumerate}
\item Every indecomposable summand of a product of special neg-isolated modules in $\Ccal$ is special neg-isolated in $\Ccal$.
\item $S$ is a right coherent ring. 
\item $C_0$ is an elementary cogenerator.
\end{enumerate}
\end{theorem}
\begin{proof}
(1) Suppose   $M$ is an indecomposable module in $\Prod(\Mcal_C)$ which is not special. Recall that $\Ycal={}^{\perp_{0,1}} C_1$ is the essential image of $\lambda_*$. By \cite[Proposition 4.15]{angeleri:hrbek:2019}  
there is an exact sequence $$0\to M\stackrel{\eta}{\rightarrow} M\otimes_R S\to M''\to 0$$ where  $M''$ belongs to $\Qcal$ and
 $\Hom_R(\eta, Y)$ is an isomorphism for every module $Y\in\Ycal$. 
In particular,  
$\Hom_R(\eta, C_0)$ is an isomorphism. Furthermore, $\Ext^1_R(M\otimes_R S, C_0)=0$ because $M\otimes_R S\in\Ycal\subseteq{}^{\perp_{1}} C_1={}^{\perp_{1}} C$  and  $C_0\in\Prod(C)$. It follows that $\Ext^1_R(M'', C_0)=0$, and of course also $\Hom_R(M'', C_0)=0$. But then 
$M''$ lies in the left perpendicular category of $C_0$, and  by Proposition~\ref{colimit of wide} it is a direct limit of modules in $\Wcal$, the class of all modules with a finite filtration by torsion, almost torsionfree modules. 

Consider now $W\in\Wcal$. We have $\Hom_R(W, M\otimes_R S)=0$  because $W$ is torsion and $M\otimes_R S$ is torsionfree. Moreover, $\Ext^1_R(W, M)=0$, because $M$ is an indecomposable, non-special module in $\Prod(C)$,  thus becomes indecomposable injective in $\Hcal_\tau$ and satisfies $\Ext^1_R(S, M)\cong\Hom_{\Hcal_\tau}(S[-1],M)=0$ for all torsion, almost torsionfree modules $S$. We conclude  that $\Hom_R(W, M'')=0$.

In conclusion, we have shown that $M''=0$ and  $M\cong M\otimes_R S$ belongs to $\Ycal$. 
Bu then $M\in{}^{\perp_{0}} C_1$, which contradicts the assumption $M\in \Prod(\Mcal_C)=\Prod(C_1)$.

(2) By Example~\ref{locally fg}, the heart  $\Hcal_\tau$ is a locally coherent Grothendieck category.  Then we know from \cite{herzog:1997},\cite{krause:1997} that the (isoclasses of) indecomposable injective objects form a topological space  $\Spec{\Hcal_\tau}$, where a basis of open subsets is given by the collection $${\mathcal O}(M)=\{E\in\Spec{\Hcal_\tau}\mid \Hom_{\Hcal_\tau}(M,E)\not=0\},\; M\in\Hcal_\tau\cap\Dcal^b(\mod R).$$
Moreover, there is a one-one-correspondence between the open subsets of $\Spec{\Hcal_\tau}$ and the hereditary torsion pairs of finite type in $\Hcal_\tau$, which maps a  torsion pair $(\Scal, \Rcal)$ to the set of indecomposable injectives  which are not contained in $\Rcal$.

Let us now consider the set $\Omega=\{ S[-1]\mid S \text{ torsion, almost torsionfree} \}$ of all simple torsionfree objects in  $\Hcal_\tau$. It generates a hereditary torsion pair $(\Scal,\Rcal)$ of finite type in $\Hcal_\tau$, which is associated to the open set $$\mathcal O=\{E\in\Spec\Hcal\mid \Hom_\Hcal(M,E)\not=0 \text{ for some } M\in \Omega\}=\Mcal_C$$ in  $\Spec{\Hcal_\tau}$.
Moreover, $\Scal$ is a localizing subcategory of $\Hcal_\tau$, and the corresponding quotient category $\Hcal_\tau/\Scal$ is again a locally coherent Grothendieck category whose spectrum is formed by the indecomposable injective objects in $\Rcal$, that is, by the complement $\mathcal O^c$ of $\mathcal O$, cf.~\cite[Thm.~2.16 and Prop.~3.6]{herzog:1997}.

We have seen in (1) that an indecomposable module in $\Prod(C)$ that is not in $\Mcal_C$ can't  belong to $\Prod(\Mcal_C)$ and thus must lie in $\Prod(C_0)$ by Proposition~\ref{no common}. That is, every indecomposable injective in $\Rcal$ is contained in $\Prod(C_0)$.  Since every injective in $\Rcal$ is a direct summand of a product of indecomposables, it follows that the injectives in $\Rcal$ coincide with $\Prod(C_0)$. In other words, the class of injective objects of $\Hcal_\tau/\Scal$ is $\Prod(C_0)=\Prod (D(S))$. Observe that these are precisely the indecomposable injective right $S$-modules. Since Grothendieck categories   are determined up to equivalence by their injective objects (cf.~\cite[p.81]{stovicek:2014}), it follows that  
$\Hcal_\tau/\Scal$ and $\Mod S$ are equivalent, and in particular, $S$ must be right coherent.

(3) Recall that a module $M$ is said to be fp-injective if every short exact sequence starting at $M$ and ending at a finitely presented module is split exact, or equivalently, if $M$ is a pure submodule of an injective module.  Since $S$ is right coherent, we know from \cite{stenstroem:1970} that the class of fp-injective right $S$-modules is closed under direct limits. In other words, the class $\Cogen_*(D(S))$ of all pure $S$-submodules of products of copies of the injective cogenerator $D(S)$ is a definable subcategory of $\Mod S$. Observe that the functor $\lambda_\ast:\Mod S\to \Mod R$ takes $\Cogen_*(D(S))$ to the class 
$\Cogen_*(C_0)$ of all pure $R$-submodules of  products of copies of $C_0$. Indeed, this is due to the fact that   $\lambda_\ast$ commutes with direct products and direct limits, and that pure-exact sequences are direct limits of split exact sequences. Conversely,   any module $X$ in $\Cogen_*(C_0)$ is an $S$-module, and any pure embedding of $X$ in a product of copies of $C_0$ is also a pure embedding of $S$-modules by \cite[Cor.~6.1.11]{prest:2009}, so that $X$ lies in $\Cogen_*(D(S))$. Now we can conclude  by \cite[Thm.~6.1.11]{prest:2009} that $\Cogen_*(C_0)$ is a definable subcategory of $\Mod R$, that is, $C_0$ is an elementary cogenerator.
\end{proof}

\begin{remark} Let $C$ be a minimal cotilting module with associated ring epimorphism $\lambda:R\to S$.
Then $S$ is right noetherian if and only if $C_0$ is $\Sigma$-pure-injective. Indeed, the only-if-part follows from the fact that $C_0$ is an injective cogenerator of $\Mod S$ and $\Prod(C_0)$ is thus closed under coproducts. Assume conversely that $C_0$ is $\Sigma$-pure-injective. By \cite{gruson:jensen:1976}  there is a cardinal $\kappa$ such that  $\Prod (C_0)\subseteq \Add (M)$ where $M$ is the direct sum of a set of representatives of  modules in $\Prod (C_0)$ of cardinality at most $\kappa$. It is not difficult to see that  $\Prod (C_0)$ then even equals $\Add (M)$ and is therefore  closed under coproducts, which shows that $S$ is right  noetherian. \end{remark}

\begin{example}\label{S not noetherian}
{\bf (1)} If $R$ is a commutative noetherian ring, then we know from Proposition~\ref{commnoeth} that $C_0$ is injective and thus $\Sigma$-(pure-)injective.

\smallskip

{\bf (2)} Consider the hereditary noetherian ring $R=\left(\begin{array}{cc} \mathbb Z&\mathbb Z\\0&\mathbb Z\end{array}\right)$. It is shown in \cite[Example 4]{beachy:1981} that the non-noetherian ring $S=\left(\begin{array}{cc} \mathbb Z_{(p)}&\mathbb Q\\0&\mathbb Z_{(q)}\end{array}\right)$, where $p,q$ are two prime numbers,  is a universal localization of $R$ at a set of matrices $\Sigma$. Hence for the minimal cotilting module $C$ given by the pseudoflat injective ring epimorphism $R\to S$ we have that $C_0$ is not $\Sigma$-pure-injective.
\end{example}
\bibliographystyle{plain}
\bibliography{bib_simples}

\end{document}